\documentclass[11pt]{amsart}
\usepackage{graphicx, amssymb,ucs}
\numberwithin{equation}{section}
\usepackage{float}
 \usepackage[utf8x]{inputenc}
 \usepackage{graphics}
\usepackage{color,ulem}

\font\tenms=msbm10
\font\sevenms=msbm7
\font\fivems=msbm5
\newfam\bbfam \textfont\bbfam=\tenms \scriptfont\bbfam=\sevenms
\scriptscriptfont\bbfam=\fivems

\newtheorem{Thm}{Theorem}
\newtheorem{Def}{Definition}[section]
\newtheorem{Rem}[Def]{Remark}
\newtheorem{Prop}[Def]{Proposition}

\newtheorem{Lem}[Def]{Lemma}
\newtheorem{Rmk}[Def]{Remark}

\def\d {{\partial}}

\newcommand{\ba}{\begin{aligned}}
\newcommand{\ea}{\end{aligned}}
\newcommand{\be}{\begin{equation}}
\newcommand{\ee}{\end{equation}}

\def \N{{\mathbf N}}
\def \R{{\mathbf R}}

\def \eps{{\varepsilon}}
\def \e{{\varepsilon}}

\def \AA{{\mathcal A}}

\def \PP{{\mathcal P}}

\def \d{{\partial}}

\begin{document}

\title[Semiclassical and spectral analysis of oceanic waves]{Semiclassical and spectral analysis of oceanic waves}

\author[Ch. Cheverry]{Christophe Cheverry}

\address[Ch. Cheverry]{Institut Math\'ematique de Rennes, 
Campus de Beaulieu, 263 avenue du G\'en\'eral Leclerc CS 74205 
35042 Rennes Cedex\\FRANCE}
\email{christophe.cheverry@univ-rennes1.fr}
\author[I. Gallagher]{Isabelle Gallagher}
\address[I. Gallagher]%
{Institut de Math{\'e}matiques UMR 7586 \\
      Universit{\'e} Paris VII \\
175, rue du Chevaleret\\
75013 Paris\\FRANCE}
\email{Isabelle.Gallagher@math.jussieu.fr}
\author[T. Paul]{Thierry Paul}
\address[T. Paul]{CNRS and Centre de Math\'ematiques Laurent Schwartz,
\'Ecole Polytechnique, 91128 Palaiseau Cedex, France}
\email{paul@math.polytechnique.fr}
\author[L. Saint-Raymond]{Laure Saint-Raymond}
\address[L. Saint-Raymond]%
{Universit\'e Paris VI and DMA \'Ecole Normale Sup\'erieure, 45 rue d'Ulm, 75230
Paris Cedex 05\\FRANCE }
\email{Laure.Saint-Raymond@ens.fr}

 \begin{abstract}
 In this work we prove that the shallow water flow, subject to strong wind forcing and linearized around an adequate stationary profile, develops for large times closed trajectories due to the propagation of Rossby waves, while   Poincar\'e waves are shown to disperse.  
 The methods used in this paper involve semi-classical analysis and dynamical systems for the study of Rossby waves, while some refined spectral analysis is required for the study of Poincar\'e waves, due to the large time scale involved which is of diffractive type.
\end{abstract}

\maketitle

\section{Introduction}

The problem we consider  is motivated by large-scale oceanography: the main physical phenomenon leading this study is    
the existence of persistent oceanic eddies, which are coherent structures of vortex type, spreading over dozens of kilometers 
and propagating slowly over periods from one year to one decade. These structures have been observed long past by 
physicists~\cite{gill,gill-longuet,Greenspan,P1,P2} who gave heuristic arguments (reproduced below) to explain their 
formation due both to   wind forcing and to   convection by a macroscopic zonal current. Giving a (much less precise) mathematical 
counterpart of those arguments, even at a linear level, requires   careful multiscale analysis and rather sophisticated tools of 
semiclassical and microlocal analysis. In this paper we simplify the model by considering particular macroscopic currents, 
which are stationary solutions of the forced equations. This allows to exhibit trapped Rossby waves, by solving the dynamics 
associated with an appropriate integrable Hamiltonian system. We prove also that the other waves produced by the dynamics, 
namely Poincar\'e waves, disperse on the same time scales (which turn out to be of diffractive type).

 \subsection{Physical observations}
  Simple observations show that {\bf large-scale ocean dynamics} can be decomposed  as the sum of   the solid-body rotation 
  together with the Earth,   convection by  macroscopic currents  (such as the Gulf Stream in the North Atlantic, the Kuroshio in 
  the North Pacific, Equatorial or Circumpolar currents), and  motion on smaller geographical zones, due for instance to the 
  fluctuations of the wind  and more generally to the coupling with the atmosphere. While the spatial extent of macroscopic 
  currents is of the order of a hundred to a thousand kilometers, those fluctuations are typically on dozens of kilometers.
We therefore expect eddies to be particular forms of those fluctuations. The point is to understand why they are quasi-stationary, 
or in other words why they do not disperse as other waves. At this stage we have to describe briefly the {\bf different kinds of 
waves} that can be found in the ocean as linear responses to exterior forcing. They are usually classified into two families, 
depending on their typical period and on their dynamical structure. The exact dispersion relation of all these waves can 
be computed explicitly \cite{cdgg,DU,GSR1,GSR2,majda} in simplified cases (no convection, linear approximation of the Coriolis 
parameter).
 \begin{itemize}
 \item {\bf Poincar\'e waves}, the period of which is of the order of a day,  are fast dispersive waves. They are due to the Coriolis 
 force, that is to the rotation of the Earth ;
 \item {\bf Rossby waves} propagate much slower, since the departure from geostrophy (that is equilibrium between pressure and Coriolis force) is very small. They are actually related to the variations of the Coriolis parameter with latitude. In particular, they propagate only eastwards.
 \end{itemize}
The  heuristic argument leading to the existence of quasi-stationary coherent structures  is then as follows (as suggested by physicists):  the wind forcing produces waves, in particular Rossby waves which would propagate, in the absence of convection, with a speed comparable to the bulk velocity of the fluid~$\bar v \sim 10 \, {\rm ms}^{-1}$; the convection by zonal flow may then stop the propagation, creating ventilation zones which are not influenced by external dynamics, in particular by continental recirculation.
We are then led to studying wave propagation under the coupled effects of the pressure, the Coriolis force and zonal convection, that is  to studying  a system of linear PDEs with non constant coefficients.

\subsection{The model}\label{themodel} The system we will consider is actually a toy model insofar as many physical phenomena are neglected. Our aim here is only to get a qualitative mechanism to explain the trapping of Rossby waves.
 More precisely,  we consider the ocean as an
{\bf incompressible, inviscid fluid with free surface} submitted to gravitation and wind forcing,
and further make the following classical assumptions~:  the density of the fluid is homogeneous~$\rho =\rho_0 =\hbox{constant } $ ; the pressure law is given by the hydrostatic
approximation $ p=\rho_0  g z$ ;  the motion is essentially horizontal and does not
depend on the vertical coordinate, 
leading to the so-called
{\bf  shallow water approximation}.
 For the sake of simplicity, we shall not discuss the effects of the interaction with the boundaries, describing neither the vertical boundary layers, known as Ekman layers, nor the lateral boundary layers, known as { Munk and Stommel layers}. We    consider a purely horizontal model, and assume an infinite domain  for the longitude (omitting the stopping conditions on the continents) as well as  for the latitude (this may be heuristically justified using the exponential decay of the equatorial waves to neglect the boundary).
The evolution of the water height $h$ and velocity $v$ is then governed by
the {\bf Saint-Venant equations with Coriolis force}
\begin{equation}
\label{SW}
\begin{aligned}
\d_t( \rho_0 h) +\nabla\cdot (\rho_0hv) =0\\
\d_t (\rho_0hv) +\nabla\cdot  (\rho_0 hv \otimes v) + \omega(\rho_0 hv)^\perp + \rho_0 g h\nabla h =\rho_0h\tau
\end{aligned}
\end{equation}
where $\omega $ denotes the vertical component of the Earth rotation vector $\Omega$, $v^\perp := (-v_2,v_1)$, $g$ is the gravity
 and
$\tau$ is the - stationary - forcing responsible for  the macroscopic flow. It depends  in particular on time averages of the wind forcing, temperature gradients and topography. The equations are written in cartesian coordinates~$(x_1,x_2)  $, where~$x_1$ corresponds to the longitude, and~$x_2$ to the latitude (both will be chosen in~$\R $). The vertical component of the Earth rotation is therefore~$\Omega \sin (x_2/R)$, where~$R$ is the radius of the Earth, but it is classical in the physical literature
to consider the linearization of~$\omega $ (known as the betaplane approximation)~$\omega(x_2) =\Omega  x_2/R$; most of our results will actually hold for more general functions~$\omega $, but in some situations we shall particularize the betaplane case in order to improve on the results.
In order to analyze the {\bf influence of the macroscopic convection} on the trapping of Rossby waves, we will  consider  small fluctuations  around the stationary solution
$$ h=\bar h ,\quad \nabla\cdot  (\bar v\otimes \bar v) +  {\omega} \bar v^\perp  =\tau  , \quad \mbox{div} \: \bar v = 0,$$
where~$\bar h$ is  a constant. 
Physical observations show  that the nonlinear convection term is essentially negligible compared to the 
Coriolis term, so that the previous equation is nothing else than the Sverdrup relation (see~\cite{P2}). 

\subsection{Orders of magnitude and scaling}
   Let us introduce the observation  length,  time and velocity scales $l_0$ (of the size of the radius of the Earth~$R$), 
   $t_0$ and $v_0$, and the nondimensional variables $\tilde x = x/l_0,$ $\tilde t = t/t_0,$ and $ u = (v-\bar v)/v_0\,.$
We also define the typical height variation $\delta h$ and the corresponding dimensionless variable
$\displaystyle \eta =  (h-\bar h) / \delta h $. We denote by $v_c$ the typical value of the velocity of  the macroscopic current: $
\displaystyle \bar  u = (\bar  u_1,0) =\bar v /v_c\,.$
The length scale $l_0$ and the convection velocity $v_c$ are fixed by the macroscopic flow: typical values for the Gulf Stream 
are~$l_0\sim 10^4 \, \hbox{km}$ and $v_c \sim 10 \,  \hbox{ms}^{-1}\,.$
As we are interested in structures  persisting during many months, 
a relevant choice for the observation  time scale is~$ t_0 =10^6 \,  \hbox{s} \: (\sim 0,38  \,  \hbox{months})\,.$
The  associated {\bf Rossby number} is then~$\displaystyle \mbox{Ro}:= 1/  (t_0 |\Omega|) =0.01$, 
recalling that~$ | \Omega | = 7.3 \times 10^{-5} \: {\rm s}^{-1}$.
The variations of water height which can be observed are typically of the order~$\delta h \sim 1 \, \hbox{m}$   
{ to be compared to } $\bar h \sim 10^3 \,  \hbox{m} \,.$
The influence of gravity (through hydrostatic pressure) is measured by  the {\bf Froude number}~$ \displaystyle\mbox{Fr}^2 := (v_0 l_0) 
/ (t_0g\delta  h) \sim 0.1, $
considering namely fluctuations of order~$v_0\sim 0.1\hbox{ms} ^{-1}.$
Defining~$\e:= \mbox{Fr}^2$ and dropping the tildas (note that as often in Physics, $\e$ is not really a very small number), we therefore end up with the following scaled system
\begin{equation}
\label{SV2}
\begin{aligned}
  \d_t \eta +  \frac1{\eps}  \nabla \cdot u +  \bar u \cdot \nabla \eta +\eps^2 \nabla \cdot (\eta u)=0\,,\\
 \d_t u  +      \frac1{\eps^2}  b  u^\perp  +\frac1\eps  \nabla \eta +\bar u \cdot \nabla u + u\cdot \nabla \bar u +\eps^2 u \cdot \nabla u =0\,,
\end{aligned}
\end{equation}
where~$b:=\omega/|\Omega|$.
 We shall compute the response to the wind forcing,  assuming  that the wind induces a {\bf pulse at  time 
$ t=0 $}: since the wind undergoes oscillations  on small spatial scales, the  initial data is further assumed to  depend both on $x$ and $x/\eps$. Typically
\begin{equation}
\label{initial-data}
(\eta_\eps,u_\eps)_{|t=0} =  (\eta_k(x), u_k(x))\exp \left( i{k\cdot x\over \eps}\right) \,,
\end{equation}
for some~$k \in{\mathbb Z}^2$. More generally we shall consider initial data which are microlocalized (in the sense of Appendix~\ref{semiclassic}) in  some compact set of~$T^* \R^{2}   $.
  
    \subsection{Local well-posedness}
   The local existence of a solution to the scaled Saint-Venant  Coriolis system (\ref{SV2}) supplemented with initial data in the form (\ref{initial-data})
   comes from the general theory of hyperbolic quasilinear symmetrizable systems.
 Defining the sound speed~$u_0$ by
$$ \eta  = \bigl \lbrack (1+ \e^3 u_0 / 2 )^2 - 1 \bigr \rbrack / \e ^3 \,, $$
we indeed obtain that (\ref{SV2}) is equivalent to
\begin{equation}\label{syst}
\displaystyle \eps^2 \partial_tU + A(x_2,\eps D)U+   \eps^3 Q(U) =  0 \, , \qquad 
U = (u_0,u_1,u_2)
\end{equation}
where  $A(x_2,\eps D_x) $ is  the linear propagator 
\begin{equation}\label{defA}
A (x_2,\eps D):=   \left( \begin{matrix}   \eps\bar u\cdot \eps \nabla  &\eps \d_1&\eps \d_2 \\
\eps \d_1 & \eps\bar u\cdot \eps \nabla  & -b( x_2)+\eps^2 \bar u_1'\\
\eps \d_2 & b( x_2) &  \eps\bar u\cdot \eps \nabla  \end{matrix} \right) \,,
\end{equation}
and  $ Q(U) := S_1 (U) \e \d_1 U + S_2 (U) \e \d_2 U $ with
\begin{equation}\label{defB}
\ S_1 (U) := \left( \begin{array}{ccc}
u_1 & \frac{1}{2} u_0 & 0 \\  
 \frac{1}{2} u_0 & u_1 & 0 \\
 0 & 0 & u_1  \end{array} \right)  \hbox{ and } S_2 (U)  := \left( \begin{array}{ccc}
u_2 & 0 & \frac{1}{2} u_0 \\  
0 & u_2 & 0 \\
 \frac{1}{2}u_0 & 0 & u_2 \end{array} \right) . 
 \end{equation}
Because of the specific form of the initial data, involving fast oscillations with respect to $x$,
we introduce semi-classical Sobolev spaces
$$
H_\eps^{s } =\{ U\in L^2\,/\, \| U\|_{H_\eps^{s}}<+\infty\} \hbox{ with } \| U\|^2_{H_\eps^{s}}=\sum_{|k|\leq s} \| (\eps \nabla)^k U\|_{L^2}^2\,.$$
We shall also need in the following to define weighted semi-classical Sobolev spaces (in the spirit of~\cite{DU}), adapted to the linear propagator as explained in Section~\ref{nonlinear}:
\begin{equation}\label{defhs1s2}
W_\eps^{s } := \Bigl\{f\in L^2(\R^2) \,/\, (1-\e^2 \partial_1^2)^{\frac{s }2 } (1-\e^2 \partial_2^2+ b^2(x_2))^{\frac{s }2 } f  \in L^2(\R^2) \Bigr\}.
\end{equation}
A classical result based on the Sobolev embedding (see Section~\ref{nonlinear} for related results)
$$\| \e \nabla U\|_{L^\infty} \leq \frac C\eps  \| \nabla U\|_{H^s_\eps} \hbox{ for any } s>1\,,$$
  implies that (\ref{SV2}) has a unique local solution $U_\eps  \in L^\infty([0,T_\eps), H^{s+1}_\eps)$.
Note that the life span of $U_\eps$ depends a priori on $\eps$. One of the goals of this article is to show     existence on an~$\e$-independent time interval.

\section{Main results and strategy of the proofs}
 Most of this paper is concerned with the analysis of the   solution to the linear equation
     \begin{equation}\label{linsyst}
\displaystyle \eps^2 \partial_tV + A(x_2,\eps D)V =  0 \, , \qquad 
V = (v_0,v_1,v_2)
\end{equation} 
which is expected to dominate the dynamics since we consider  small fluctuations. The description of the linear dynamics is provided in Theorem~\ref{mainresult} below.
The comparison between linear and nonlinear solutions is postponed to the final section of the paper (see Theorem~\ref{nonlineartheorem}).

For technical reasons we shall restrict our attention in this paper to the case of a shear flow, in the sense that~$
\bar u (x)= (\bar u_1(x_2),0),
$
where~$\bar u_1$ is a smooth, compactly supported function. We shall further assume for simplicity that the zeros 
of~$\bar u_1$,  in the interior of its support, are {\bf of order one}. We shall also suppose throughout the paper that~$b$ is a smooth function with a  symbol-like behaviour:
\begin{equation}\label{bsymbol}
\forall \alpha \in \N,  \:\exists C_\alpha, \: \forall y \in \R, \quad |b^{(\alpha)} (y)| \leq C_\alpha (1+b^2(y)),
\end{equation}
and we shall further assume that
$
\displaystyle \lim_{y \to \infty} b^2(y) = \infty,
$ 
and that~$b$   has at most a finite number of   critical points (that is to say points where~$b'$ vanishes).  Without such an assumption one could not  construct Rossby waves.
We shall also  suppose that the initial data is {\bf microlocalized} (see Appendix~\ref{semiclassic}) in  some compact set of~$T^* \R^{2}   $  (which we shall identify to~$\R^4$ in the following), denoted~$\mathcal C$  and satisfying
\begin{equation}\label{cond1}
\mathcal C\cap\{\xi_1=0\} = \emptyset .
\end{equation}
Thanks to this assumption, which is propagated by the linear flow, one can diagonalize the system into Rossby and Poincar\'e modes.
Finally in order to avoid pathological  trapped  Rossby trajectories we shall
also require that
 \begin{equation}\label{cond3}
\mathcal C\cap \Sigma= \emptyset ,  
\end{equation}
where~$\displaystyle \Sigma$ is a codimension 1 subset of~$\R^4$  defined in Proposition~\ref{classification}.
    \subsection{Statement of the main results}
    In this paragraph we shall state the two main theorems proved in   this paper.
   The first result deals with the linear system~(\ref{linsyst}).
    \begin{Thm}[The linear case] \label{mainresult} 
 There is  a submanifold~$\Lambda$  of~$\R^4$, invariant under  translations in the $ x_1 -$direction, such that the following properties hold.

  Let $U_{\e,0} $ be $ \e -$microlocalized in
a compact set~${\mathcal C}$ satisfying Assumptions~(\ref{cond1})-(\ref{cond3}). For any parameter~$\eps >0$, denote by $V_\eps$ the associated solution to (\ref{linsyst}).
Then for all~$t \geq 0$ one can write~$V_\e(t)$ as the sum of a ``Rossby" vector field and a ``Poincar\'e" vector field: $ V_\e(t)  = V_\e^R(t)  + V_\e^P(t) $, satisfying the following properties:
\begin{enumerate}
\item \label{rossbybehaviour}
There is a compact set~$K $ of~$\R^2$  such that
$$\forall t\geq 0, \quad  \|    V_\e^R(t)\|_{L^2(K)}\neq O(\e^\infty) 
$$
if and only if the~$\e$-frequency set of~$ V_\e^R(0) $     intersects~$\Lambda  $.
\item \label{poincarebehaviour}
Suppose   that~$b^2$ has only one non degenerate critical value (meaning that~$(b^2)'$ only vanishes at one point, where~$(b^2)''$ does not vanish). Then  for any compact set~$\Omega$ in~$\R^2$, one has
$$
\forall  t> 0,\quad  \| V_\e^P(t)\|_{L^2(\Omega)}=O(\e^\infty).
$$
\end{enumerate}
In particular 
 supposing  that~$b^2$ has only one non degenerate critical value, then  there is a compact set~$K $ of~$\R^2$  such that
$$\forall t>0, \quad  \|  V_\e(t)\|_{L^2(K)}\neq O(\e^\infty) 
$$
if and only if the~$\e$-frequency set of~$ V_\e^R(0)$     intersects~$\Lambda  $.  \end{Thm}
    \begin{Rem}   
 Actually $\Lambda$
 corresponds to the set of initial positions and frequencies in the phase space giving rise to trapped trajectories for the Rossby hamiltonian. This will be made more precise in Section~\ref{rossby}, where we shall prove that under some additional (non restrictive) assumptions on~$\bar u$, $\Lambda$ is of codimension one.   In particular it will be shown that some of those trapped trajectories actually exhibit a singular behaviour in large times, in the sense that they converge in physical space towards a point, while the~$\xi_2$ frequency goes to infinity. This could be interpreted like the creation of some sort of oceanic eddies.

  The result~(\ref{poincarebehaviour}) is related to   dispersive
 properties of the Poincar\'e hamiltonian on diffractive type times (of the type~$O(1/\e^2)$), which requires some spectral analysis. Due to the assumption on~$b^2$ one can write a rather simple proof; more general conditions could be treated, but the Bohr Sommerfeld quantization would require to decompose the phase space into various zones according to the geometry of the level sets of the Hamiltonian, which is much more technical and beyond the scope of this article.  Actually in~\cite{gpsrmourre} we propose a different approach, based on Mourre estimates, which allows to relax very much the assumptions on~$b^2$ and on~$\bar u$.
  \end{Rem}
 The final section of this paper is devoted to the proof of the following theorem, which states that the very weak coupling chosen in this paper implies that nonlinear dynamics are governed by the linear equation.    We also consider more generally the following weakly nonlinear system (with the notation~(\ref{defA}) and~(\ref{defB})):
   \begin{equation}\label{WNLS}
\eps^2 \d_t U+ A(x,\eps D_x) U +\eps ^{3+\eta} S_1(U) \eps \d_1 U +\eps^{3+\eta} S_2(U) \eps \d_2 U=0 ,\quad \eta\geq 0\,.
\end{equation}
The case~$\eta = 0$ corresponds of course to the original system~(\ref{syst}) presented in the introduction.
  \begin{Thm}[The nonlinear case] \label{nonlineartheorem} 
    Let $U_{\e,0}$ be any initial data bounded in $W_\eps ^{4}$. Then the following results hold.
  \begin{enumerate}
  \item The case~$\eta = 0$:
    \begin{enumerate}
    \item \label{uniformlifespan} There exists some $T^*>0$ such that
the initial value problem~(\ref{WNLS}) with~$\eta = 0$ 
 has a unique solution~$ U_\eps$ on $[0,T^*[$ for any $\eps >0$.
  \item\label{iftheorem} Assume that the solution $V_\eps$ to the linear equation (\ref{linsyst}) satisfies
 $$\| \eps V_\eps \|_{L^2([0,T^*[;L^\infty)} \to 0 \hbox{ as } \eps \to 0.$$
 Then the solution~$U_\eps$ to~(\ref{WNLS}) with~$\eta = 0$ satisfies
 $$\| U_\eps-V_\eps \| _{L^2} \to 0 \hbox{ uniformly  on }[0,T^*[ \hbox{ as }\eps \to 0.$$
\end{enumerate}
 \item \label{weakercoupling}  The case~$\eta>0$: 
 Le~$T>0$ be fixed. Then there is~$\e_0>0$ such that for any~$\e \leq \e_0$, the equation~(\ref{WNLS})  has a unique solution~$U_\eps$ on~$[0,T]$. Moreover,
  $$\| U_\eps-V_\eps \| _{L^2} \to 0 \hbox{ uniformly  on }[0,T]\hbox{ as }\eps \to 0.$$
  \end{enumerate}
  \end{Thm}
\begin{Rem}
 Result~(\ref{iftheorem}), joint with Theorem~\ref{mainresult}, implies in particular that as soon as~$b^2$ has only one non degenerate critical value, then for positive times the energy of~$U_\eps$ on any fixed compact subset is carried only by Rossby waves. The refined $L^\infty$ estimate on the linear solution required in  result~(\ref{iftheorem}) should be proved by using WKB tools.
For the sake of simplicity, we shall not consider such technical estimates here, all the less that we do not expect them to be enough to get an optimal result regarding the nonlinear problem (see Remark \ref{rk-embedding}). 
That is the reason why we consider, in result~(\ref{weakercoupling}), a weaker coupling still. That result implies in particular that  the $L^2$ norm of $U_\eps$ on any fixed compact subset may remain bounded from below only if there are trapped Rossby waves, i.e. only  if the $\eps$-frequency set of the initial data does intersect $\Lambda$ (with the notation of Theorem~\ref{mainresult}).
\end{Rem}
   \subsection{Some related studies}
   This work follows a long tradition of mathematical studies of fast rotating fluids, following~\cite{schochet} and~\cite{grenier}; we refer for instance to~\cite{cdgg} and~\cite{GSR2}   for a number of references. 
    The present study  concerns the case when the penalization matrix does not have constant coefficients. A first study in this type of situation may be found in~\cite{ftvariablegsr}, where a rather general penalization matrix was considered. Due to the generality of the situation, explicit computations were ruled out  and no study of waves was carried out. In order to     compute    explicitly the modes created by the penalization matrix, various authors (see~\cite{DM}, \cite{DM2}, \cite{DU} as well as~\cite{GSR1})  studied the betaplane approximation, in which the rotation vector depends linearly on the latitude. In that case explicit calculations may again be carried out (or some explicit commuting vector fields may be computed) and hence again one may  derive envelope equations. In this paper we choose again to work with a more general rotation vector, this choice being made possible by a semi-classical setting (see the next paragraph); in particular that setting enables us diagonalize the system approximately, therefore to compute waves; note that a related study is performed by two of the authors in~\cite{CP1} via a purely geometric  optics approach, where no explicit diagonalization is performed (actually  the initial data is strongly polarized so that only Rossby modes are present, including in the non linear setting) .
   
Another feature of our study is that it is a multi-scale problem, in the sense that the oscillation frequency is much bigger than the variation of the coefficients of the system. This is dealt with by using semi-classical analysis (which, compared to the previous paragraph, enables us to compute {\it almost commuting} vector fields although the penalization matrix no longer depends only linearly on the latitude).
      Such techniques are classical in geometrical optics, but in our case the additional difficulty is that the propagators  are linked to different time scales: this is due to the fact that the system has eigenvalues at different scales (one is actually a subsymbol). In particular   we are mostly interested in the role of the subsymbol in the dynamics, as this subsymbol is responsible for the trapping phenomenon we want to exhibit: this implies, by  semi-classical analysis, the need to study the dynamical system induced by that subsymbol. On the other hand this also means that the dynamics linked to  other eigenvalues 
     must be analyzed on diffractive-type time scales, therefore much longer than that allowed by semi-classical analysis. We are able to show the dispersion of those waves by using spectral analysis and Bohr-Sommerfeld quantization.

   \subsection{Organization of the paper}
The proof of Theorems~\ref{mainresult} and~\ref{nonlineartheorem} requires a number of steps which are described in this paragraph.
 \subsubsection{Reduction to scalar propagators}
 Persistent structures are related to the propagation of Rossby waves. Our first task is therefore to transform the original linear system~(\ref{linsyst}) into three scalar equations. One is polarized on Rossby waves while the two others are polarized on  Poincar\'e waves. This is done in Section~\ref{scalar} by proving some necessary conditions
 for the existence of those propagators. The general strategy is the following:
 \begin{enumerate}
\item Consider the system~$  A(x_2,\e D)U =i \tau U$. Take the Fourier transform in~$x_1$, which is possible since the equation is translation invariant in~$x_1$. Then extract from this system  (by linear combinations and substitutions) a linear equation on one component~$u_k$ of~$U$,  of the type~$h(x_2,\e D_2;\xi_1, \e, \tau)u_k = 0$.
\item The symbolic equation corresponding to the PDE writes $h(x_2,\xi_2;\xi_1, \e, \tau) = 0$. It has  three roots
(with respect to $ \tau $) , $ \tau_\pm(x_2,\xi;\eps)$ (Poincar\'e roots) and~$\tau_{   R }(x_2,\xi;\eps)$ (Rossby root). We find
$ \tau_\pm (x_2,\xi;\eps) = \tau_\pm (x_2,\xi) + O(\eps) $ and $\tau_{ R} (x_2,\xi;\eps) = \eps \, \widetilde \tau_R 
(x_2,\xi)+ O(\eps^2)$. 
\item \label{step4} Those roots are not necessarily symbols. 
To guarantee these are indeed symbols one needs a microlocalization. Given a compact set and a truncation on that compact set~$\chi$, one can construct three operators~$T_j^{\chi}$ (via a  general theorem, stated and proved in an abstract way in Theorem~\ref{lemkiditou} of Appendix~\ref{sectionlemkiditou}) whose principal symbols are 
  precisely the Rossby and Poincar\'e symbols~$ \tau_\pm$ and~$ \tau_{ R}$.  \end{enumerate}

\subsubsection{Trapping of Rossby waves}\label{introrossby}
Section~\ref{rossby} is devoted to the study of the Rossby propagator~$T_0$, and in particular to the proof of result~(\ref{rossbybehaviour})
 in Theorem~\ref{mainresult}. For the time scale   considered, it is easy to see that the energy propagates according to the trajectories of the semiclassical Rossby hamiltonian~$\widetilde \tau_R$.  The first step of the analysis therefore consists in studying the dynamical system giving rise to those trajectories. It turns out that the  trajectories are always bounded in the~$x_2$ direction. One is therefore reduced to studying the trajectories in the~$x_1$ variable and in identifying the set~$\Lambda$ of initial data in the cotangent space giving rise to trapping in~$x_1$. 
One then checks that~$\Lambda$ is of codimension one under some additional assumptions on~$\bar u$, and the last step of the study consists in studying more precisely the trajectories in some specific situations, in particular in the case of the betaplane approximation.

  \subsubsection{Dispersion of Poincaré waves}
 The next step of our analysis of wave propagation consists in proving, in Section~\ref{sectionpoincare}, 
  that Poincar\'e waves propagate so fast that they exit from any (bounded) domain of observation on the time scale that we consider, which proves result~(\ref{poincarebehaviour})
 in Theorem~\ref{mainresult}.   Note that, because of the very long time scaling, usual tools of semiclassical analysis cannot be applied for the Poincar\'e waves~: we actually need deeper arguments such as the Bohr-Sommerfeld quantization to conclude.

\subsubsection{A diagonalization result} 
Once the Rossby and Poincar\'e propagators have been well understood, we can retrace the steps followed in Section~\ref{scalar} to prove that the necessary conditions on the scalar propagators are sufficient.  The difficulty is that the operators~$\Pi$ enabling one to go from the original system to the scalar equations and back (computing an approximate left inverse~$Q$ of~$\Pi$ at the order~$O(\e^\infty)$) are only continuous on microlocalized functions; moreover the scalar propagators~$T_j^{\chi}$  are themselves only defined on 
microlocalized functions. 
 So that requires understanding the persistence of the  microlocalization of the solutions to the scalar equations. That is achieved in the two previous sections, where it is proved that if the initial data is conveniently microlocalized, then for any time~$t\geq 0$ one can find a compact set~$K$ and one can construct~$T_j^{\chi}$ as in Section~\ref{scalar}, so that  the solution to the scalar equations with propagators~$T_j^{\chi}$ is microlocalized in~$K$ (actually for Poincar\'e modes the microlocalization is in the variables~$(x_2,\xi_1,\xi_2)$ only, which is enough for our purpose). This enables us in Section~\ref{diago} to conclude rather easily by computing explicitly the matrix principal symbols of~$\Pi$ and~$Q$.

\subsubsection{The analysis of the nonlinear equation}  Section~\ref{nonlinear} is devoted to the proof that the solution of the nonlinear equation remains close to that of the linear equation. The method of proof consists
first in proving the wellposedness of the nonlinear equation on a uniform time interval, by using semi-classical weighted  Sobolev type spaces, whose additional feature is to be well adapted to the penalization operator~$A(x_2,
\e D)$: one therefore  constructs a matrix-valued pseudo-differential operator 
which approximately commutes with~$ A(x_2,\e D) $. The convergence of~$U_\e - V_\e$ to zero relies on a standard~$L^2$
   energy estimate on~$U_\e - V_\e$,  and Gronwall's lemma.  

\subsubsection{Two appendixes} In Appendix~\ref{sectionlemkiditou} one can find the statement and the proof of a general theorem,  used in Section~\ref{scalar}, allowing to associate
to a linear evolution PDE a number of operators describing the dynamics of the equation; those operators are constructed by writing down the symbolic equation associated to the PDE and in quantizing the roots  of that polynomial (in the time derivative).
 Appendix~\ref{semiclassic} finally collects a number of prerequisites on microlocal and semiclassical analysis, 
 that are used throughout the paper.

\section{Reduction to scalar propagators}\label{scalar}
\subsection{Introduction} 
Let us first recall that the propagator
$$
A(x,\eps D ) = \left( \begin{matrix}  \eps\bar u\cdot \eps \nabla &\eps \d_1&\eps \d_2 \\
\eps \d_1 &\eps\bar u\cdot \eps \nabla& - b + \e^2 \bar u'_1 \\
\eps \d_2 & b  &\eps\bar u\cdot \eps \nabla \end{matrix} \right)
$$
can, in the particular case when $\bar u \equiv 0$ and $b(x_2)=\beta x_2$, be diagonalized without any error term  (in particular   for any  finite~$\eps$), using a Fourier  basis $(\exp (\frac i\eps x_1 \xi_1))_{\xi_1\in \R}$ in $x_1$ and a Hermite basis $(\psi_n^\eps(x_2))_{n\in \N}$ in~$x_2$. 
  Precisely,  the following statement is proved in~\cite{GSR1}.
\begin{Prop}[Gallagher \& Saint-Raymond,\cite{GSR1}]
{ 
For all $(\xi_1, n,j)\in \R\times \N\times 
\{-1,0,1\}$, denote by $\tau(\xi_1, n ,j)$ the three roots  (in increasing order in~$j$) of
\begin{equation}
\label{tau-hermite}
\tau^3 -(\xi_1^2+\beta\eps (2n+1)) \tau+\eps\beta \xi_1=0.
\end{equation}
Then there exists a complete
family of $L^2(\R\times \R,\R^3)$  of
pseudo-eigenvectors~$(\Psi^\eps_{\xi_1, n,j}) $ of the operator~$A(x,\eps D)$ (where~$\bar u \equiv 0$ and $b(x_2)=\beta x_2$):
\begin{equation}
\label{eigenvectors}
\forall ( \xi_1, n,j)\in \R\times \N\times \{-1,0,1\}, \quad A(x,\eps D)\Psi^\eps_{\xi_1, n,j} =i\tau(\xi_1, n,j) \Psi^\eps_{\xi_1, n,j}
\end{equation}
where  $\Psi^\eps _{\xi_1, n,j}$ can be computed in terms of the $n$-th Hermite function $\psi_n^\eps(x_2)$ and its derivatives.}
\end{Prop}
In other words, the three scalar propagators (numbered by~$j$) can be obtained from the symbolic equation  (\ref{tau-hermite}) remarking that $\beta\eps(2n+1)$ is  the quantization of the harmonic oscillator~$-\eps^2 \d_2^2+\beta^2 x_2^2$.  It is proved in~\cite{GSR1} that as~$\xi_1$ and~$n$ go to infinity
$$
\tau(\xi_1, n ,\pm) \sim \pm \sqrt{\xi_1^2 + \beta\eps(2 n+1)}, \quad \mbox{and} \quad \tau(\xi_1, n ,0) \sim   \frac{\eps\beta\xi_1 }{\xi_1^2 + \beta\eps(2 n+1)} \cdotp
$$
We are interested here in deriving  a symbolic equation similar to (\ref{tau-hermite}) for a general  zonal current $\bar u =(\bar u_1(x_2),0) $ and Coriolis parameter $b=b(x_2)$. The difficulty comes from the fact that  the propagation of waves is governed by  a matrix of differential operators with non-constant coefficients, the diagonalization of which is not a standard computation.
Of course, in the semiclassical limit $\eps \to 0$, we expect to get a good approximation of the propagation at leading order by considering the matrix of principal symbols
\be \label{matrixps}
A_0(x,\xi) := \left( \begin{matrix} 0 &i\xi_1&i\xi_2\\
i\xi_1 &0 & -b(x_2)\\
i\xi_2 & b( x_2) & 0 \end{matrix} \right) 
\ee
 and by computing 
  the scalar propagators associated to each eigenvalue
   \begin{equation}\label{deftaupm}
   \tau_\pm (x_2, \xi_1,\xi_2 ) :=\pm \sqrt{\xi_1^2 + \xi_2^2 + b^2(x_2)}
\end{equation}
and~0. The eigenvalue 0 corresponds to the Rossby modes, whereas the two~$O(1)$ eigenvalues~$\pm\sqrt{  \xi_1^2 + \xi_2^2  + b^2(x_2)}$ are the Poincar\'e modes.
 Nevertheless, this approximation is relevant only for times of order $O(\frac1\eps)$, and we are interested here in much longer times, of order~$O(\frac1{\eps^2})$. This means that we need to compute the next order of the expansion of the eigenvalue~$0$.  Once that is done, we need to quantify   these symbol eigenvalues to deduce scalar propagators.
 \subsection{The Rossby modes} \label{diagorossby}
Finding precisely the Rossby modes (up to an~$O(\e^2)$ error) requires more intricate calculations than merely diagonalizing the matrix of principal symbols~$A_0(x,\xi) $ given in~(\ref{matrixps}).
 So let $i \tau$ be an eigenvalue of the propagator, assumed to be of the form~$\tau = \e \widetilde\tau_R + O(\e^2)$. Explicit computations lead to the subsystem
$$
 \left( \begin{matrix} \eps \bar u_1 \xi_1 -\tau & \xi_1\\ \xi_1 & \eps \bar u_1\xi_1-\tau\end{matrix}\right) \left(\begin{matrix} \rho\\u_1\end{matrix}\right) = \left(\begin{matrix} i \eps \d_2 u_2\\ -i( b  - \eps^2 \bar u_1')  u_2\end{matrix}\right)  $$
 and defining
$
\displaystyle\alpha_\tau(x_2,\xi_1):=\eps \bar u_1 \xi_1 -\tau $ and $p_R(x_2,\xi_1) =-\xi_1^2 +\alpha_\tau^2 
$  to the scalar equation
 \begin{equation}
 \label{tau}
 \left(
 \e \partial_2 p_R^{-1}\bigl(i\alpha_\tau \e \partial_2 + i \xi_1 (b  - \eps^2 \bar u_1')\bigr) - b  p_R^{-1} (i\xi_1 \e \partial_2 +i\alpha_\tau  (b  - \eps^2 \bar u_1')) +i \alpha_\tau
 \right)u_2 = 0
  \end{equation}
where $\xi_1$ is the Fourier variable corresponding to $x_1/\eps$. From now on we assume that~$\xi_1$ is fixed, and is bounded away from zero (recalling Assumption~(\ref{cond1})).
Note that equation~(\ref{tau}) makes sense because we are assuming here that~$\tau = \e \widetilde \tau_R + O(\e^2)$, so~$ p_R^{-1}$ is well defined. That would not be the case for the Poincar\'e modes (where~$\tau = \pm\sqrt{ |\xi|^2 + b^2(x_2)} + O(\e)$) so we shall use another subsystem in the next paragraph to deal with the Poincar\'e operators.

 In order to derive a symbolic equation associated with the differential equation~(\ref{tau}), we shall proceed by transforming~(\ref{tau}) into a differential equation which is the
left quantization  (in the sense recalled in Appendix~\ref{semiclassic}) of a symbol, polynomial in~$\eps$ and~$\tau$.  This leads
to a differential equation of the type 
 \begin{equation}
 \label{tautau}
{\mathcal P}_2 (x_2  ;\eps,\tau)( \eps \partial_2)^2 u_2 + {\mathcal P}_1 (x_2;\eps,\tau) \eps \partial_2  u_2 + {\mathcal P}_0 (x_2 ;\eps,\tau) u_2   = 0,
 \end{equation}
where each~${\mathcal P}_j (x_2 ;\eps,\tau)$ is a smooth function in~$x_2$,   and has polynomial dependence  in~$\eps$ and in~$\tau$ (precisely of degree at most 5 in~$\eps$ and~$\tau$). This generalizes~(\ref{tau-hermite}); one can compute in particular, using the fact that~$\e \partial_2 \alpha_\tau = O(\e^2)$, that
$$
{\mathcal P}_2 (x_2  ;\eps,\tau) = ip_R \alpha_\tau + O(\e^2), \quad  {\mathcal P}_1 (x_2;\eps,\tau) = O(\e^2), \quad \mbox{and}$$
$$ {\mathcal P}_0 (x_2 ;\eps,\tau) = ip_R \left(\xi_1 \eps b' - (b^2+ \xi_1^2) \alpha_\tau\right) + O(\e^2).
$$
The differential operator appearing on the left-hand side of~(\ref{tautau})
 is the left quantization of   the following symbol:
  \begin{equation}
 \label{tauencore}
h(x_2, \xi_2; \eps,\tau) := -  {\mathcal P}_2 (x_2;\eps,\tau)\xi_2 ^2   + i  {\mathcal P}_1 (x_2;\eps,\tau) \xi_2    + {\mathcal P}_0 (x_2;\eps,\tau)   
 \end{equation}
 which belongs for each~$\tau$ to~$S_{2}(g_\tau)$ for some function~$g_\tau$ of the type
 $$
 g_\tau(x_2,  \xi_2) = (1+\tau^5) \left( 1+  b^2(x_2) + \xi_2^2\right)
 $$
 recalling that~$ \bar u_1$ and~$\bar u_1'$ are bounded from above, as well as Assumption~(\ref{bsymbol}).
 
Now we recall   that~$h(x_2 , \xi_2; \eps,\tau)  $ is a polynomial of degree 5 in~$\tau$ hence has five roots, among which~2 are actually spurious: these are of the form~$   \pm \xi_1 + O(\eps) $, and they appear   because we have multiplied the equation by a polynomial in~$\tau$ which cancels at the point~$   \pm \xi_1$ at first order in~$\eps$.
An easy computation allows to obtain the two other~$O(1)$ roots, which are precisely the Poincar\'e roots~$ \pm\sqrt{ |\xi|^2 + b^2(x_2)}+O(\e) $, and an asymptotic expansion allows also easily to derive the Rossby~$O(\e)$ root:  one finds
\begin{eqnarray}\label{deftau0}
  \tau_R (x_2, \xi_1,\xi_2; \eps)  &:= & \e \widetilde  \tau_R (x_2, \xi_1,\xi_2; \eps) + O(\eps^2)  , \quad \mbox{where} \nonumber \\
  \widetilde  \tau_R (x_2, \xi_1,\xi_2 ) &:=& \frac{b'(x_2) \xi_1}{\xi_1^2 + \xi_2^2 + b^2(x_2)} + \xi_1 \bar u_1(x_2)  .
\end{eqnarray}
Now that the root~$\tau_R$ has been computed,  our next task is to prove the existence of the Rossby propagator~$T_R = \eps \widetilde T_R$ whose principal symbol is precisely~$\eps  \widetilde \tau_R$. Actually this result is a direct consequence of Theorem~\ref{lemkiditou} stated and proved in Appendix~\ref{sectionlemkiditou}: with the notation of Theorem~\ref{lemkiditou}, one has~$\nu = 1$ and~$
\partial_\tau h_0 (x,\xi,0) = -i \xi_1^2 (\xi_1^2 + \xi_2^2 + b^2(x_2)).$
In the following statement, the time variable~$s$ is defined as~$s = t/\e^2$.
\begin{Prop}[The Rossby propagator]\label{therossbymode}
Let~$\widetilde\tau_R$ be the symbol defined in~(\ref{deftau0}). Then for any compact set~${\mathcal K}$ satisfying Assumptions~(\ref{cond1},\ref{cond3}) there exists a formally self-adjoint  pseudo-differential operator~$\widetilde T_{R}$  of principal symbol~$\widetilde\tau_R$  such that if~$\varphi_R$ is microlocalized in~${\mathcal K}$ and solves
\begin{equation}\label{eqrossby}
   \partial_s \varphi_R = i\e \widetilde T_R \varphi_R, 
\end{equation}
then
$$
{\mathcal P}_2 (x_2  ;\eps,\partial_s)( \eps \partial_2)^2  \varphi_R + {\mathcal P}_1 (x_2;\eps, \partial_s) ( \eps \partial_2)   \varphi_R + {\mathcal P}_0 (x_2 ;\eps, \partial_s)  \varphi_R   = O(\e^\infty).
$$
\end{Prop}
  \begin{Def}[The Rossby operator]\label{therossbyproj}
We shall call~$\Pi_R$  the Rossby operator defined  by
 $$
 \Pi_R :=\left(\begin{array}{c}
P_R^{-1}
\bigl(
 i (-i \eps^2 \bar u_1 \partial_1 -\e \widetilde T_R) \e \partial_2 + \e \partial_1(b-\e^2 \bar u_1')  \bigr) \\
-P_R^{-1}\bigl(- \e^2 \partial_{12} + i  (-i \eps^2 \bar u_1 \partial_1 -\e \widetilde T_R)(b-\e^2 \bar u_1')  \bigr) \\
\rm{Id}\end{array}\right) \, , $$ where~$ P_R := \e^2 \partial_1^2- ( i \eps^2 \bar u_1 \partial_1 +\e \widetilde T_R)^2.
$ 
\end{Def}
 \begin{Rem}\label{formaldiagorossby} \begin{enumerate}
 \item 
 Notice  that~$\Pi_R$ is well defined since the principal symbol of~$P_R$ is bounded   from below (see Appendix~\ref{semiclassic}).
  \item Proposition~\ref{classification} 
  shows that if~$\varphi_{R|t = 0}$ is microlocalized in a compact set~${\mathcal K}_0$ satisfying~(\ref{cond3}), then the solution to~(\ref{eqrossby}) is microlocalized for all~$t\geq 0$  in a compact set~${\mathcal K}_t$.
\item The above computations allow to formally recover the original shallow-water equation, up to~$O(\eps^\infty)$.   Indeed  retracing the steps which  enabled us above to derive equation~(\ref{tautau})   shows that if~$ \varphi^0  $ is a smooth function  conveniently microlocalized, and if~$\varphi$ solves
$$ \e \partial_t \varphi = i \widetilde T_R \varphi, \quad \varphi_{|t = 0} = \varphi ^0,
 $$
then  the vector field~$\displaystyle
   U:= \Pi_R \varphi $ 
satisfies (\ref{linsyst}) up to~$O(\e^\infty)$.  This property will be made rigorous in Section~\ref{diago}.
\end{enumerate}
  \end{Rem}
 \subsection{The Poincar\'e modes} \label{diagopoincare}
In this paragraph we shall follow the method used above in the case of Rossby modes to infer   Poincar\'e propagators~$T_\pm$ and operators~$\Pi_\pm$. Actually one cannot use precisely the same method since  the symbol~$ (\eps \bar u_1\xi_1-\tau)^2 -\xi_1^2 $ may vanish when~$\tau =   \tau_\pm + O(\e)$. So we shall instead consider the subsystem 
$$
 \left( \begin{matrix} \eps \bar u_1 \xi_1 -\tau & -i\eps\partial_2 \\- i\eps\partial_2 & \eps \bar u_1\xi_1-\tau\end{matrix}\right) \left(\begin{matrix} \rho\\u_2\end{matrix}\right) = \left(\begin{matrix}  -\xi_1 u_1\\ ib u_1\end{matrix}\right)  $$
  and   the scalar equation 
   \begin{equation}
 \label{taupoincare}
  \begin{aligned}
\Bigl( \xi_1 p_P^{-1} (- \alpha_\tau \xi_1 - \e \partial_2 b  )  + \alpha_\tau + (b-\e^2 \bar u'_1) p_P^{-1} ( \e \partial_2- \alpha_\tau b ) \Bigr)u_1 =0,
  \end{aligned}
  \end{equation}
  where as before~$ \alpha_\tau (x_2,\xi_1):= \eps \bar u_1\xi_1-\tau$, and
  $
 p_P  (x_2,\xi_1):= \e^2 \partial_2 ^2 +\alpha_\tau^2  .
  $
  We notice here that~$\alpha_\tau^{-1}$ is well defined when~$\tau = \tau_\pm$     since~$\tau_\pm$ is bounded away from zero by the assumption on~$\xi_1$, and the same goes for~$\tau_\pm^2 - \xi_2^2$ so~$p_P^{-1}$ is also well defined.  
 Then it remains to follow the steps of Paragraph~\ref{diagorossby} to obtain a new scalar PDE of the same type as~(\ref{tautau}), as well as a symbol equation of the type~(\ref{tauencore}):
  \begin{equation}
 \label{tauencorepoincare}
\tilde h(x_2, \xi_2; \eps,\tau) := -  \tilde {\mathcal P}_2 (x_2;\eps,\tau)\xi_2 ^2   + i  \tilde {\mathcal P}_1 (x_2;\eps,\tau) \xi_2    + \tilde {\mathcal P}_0 (x_2;\eps,\tau).   
 \end{equation}
Of course~$  \tau_\pm  $ are roots  of that equation up to~$O(\e)$, and~$\e \tau_R$ is a root up to~$O(\e^2)$. 
   Then the application of Theorem~\ref{lemkiditou} implies a similar result to Proposition~\ref{therossbymode}, noticing that with the notation of Theorem~\ref{lemkiditou}, $\nu = 0$ and
   $
   \partial_\tau h_0 (x,\xi,\tau_\pm) = -2 (\xi_1^2 + \xi_2^2 + b^2(x_2))
.
$
\begin{Prop}[The Poincar\'e propagator]\label{thepoincaremode}
Let~$ \tau_\pm$ be the symbol defined in~(\ref{deftaupm}). Consider a compact set~${\mathcal K}_P \subset \R^* \times T^*\R  $. Then there exists   formally self-adjoint pseudo-differential operators~$T_{\pm}$  of principal symbols~$ \tau_\pm$  such that if~$\varphi_ \pm $ solves
$
   \partial_s \varphi_ \pm = i T_\pm \varphi_\pm, 
$
and~$\varphi_ \pm$ is microlocalized in~$\R \times {\mathcal K}_P$, then
$$
\tilde {\mathcal P}_2 (x_2  ;\eps,\partial_s)( \eps \partial_2)^2  \varphi_\pm + \tilde{\mathcal P}_1 (x_2;\eps, \partial_s) ( \eps \partial_2)   \varphi_\pm+ \tilde{\mathcal P}_0 (x_2 ;\eps, \partial_s) \varphi_\pm   = O(\e^\infty).
$$
\end{Prop}
 \begin{Def}[The Poincar\'e operator]\label{thepoincareproj}
We shall call~$\Pi_ \pm $  the  Poincar\'e operator defined  by
 $$
 \Pi_ \pm :=\left(\begin{array}{c}
 P_ \pm ^{-1} (i( -i\eps^2 \bar u_1\partial_1-{T_\pm})  \e \partial_1 + \e \partial_2 b)\\
 \rm{Id} \\
  P_ \pm ^{-1} (i ( -i\eps^2 \bar u_1\partial_1-{T_\pm})   b- \e \partial_1   \e \partial_2 ) \end{array}\right) \,  , 
 \quad \mbox{where} \quad  P_ \pm :=   (  i\eps^2 \bar u_1\partial_1+{T_\pm})  ^2  +  \eps^2 \partial_2^2 .
$$
\end{Def}
 \begin{Rem}\label{formaldiagopoincare}  As in Remark~\ref{formaldiagorossby}, one sees formally that if~$ \varphi^0  $ is a smooth function conveniently microlocalized, and if~$\varphi$ solves
$$ \e^2 \partial_t \varphi = i T_\pm \varphi, \quad \varphi_{|t = 0} = \varphi ^ 0,
   $$
then  the vector field~$\displaystyle
   U:= \Pi_\pm \varphi $ 
satisfies (\ref{linsyst}) up to~$O(\e^\infty)$.  This property will be made rigorous in Section~\ref{diago}.
  \end{Rem}

 \section{Study of the Rossby waves}\label{rossby}
\subsection{The dynamical system}\label{dynamicalsystem}$ $
For the time scale considered here, the propagation of energy by Rossby waves is given by the transport 
equation (see Appendix~\ref{semiclassic})
$$\d_t f+\{ \widetilde\tau_R, f\} =0$$ 
where $ \widetilde \tau_R$ is the principal symbol of the Rossby mode computed in~(\ref{deftau0}):
\begin{equation}
\label{tau0-symbol}
\widetilde\tau_R (\xi_1,x_2,\xi_2)  = {b'(x_2) \xi_1\over \xi_1^2+\xi_2^2+b^2
(x_2)} +\bar u_1(x_2) \, \xi_1\,.
\end{equation}
As $\widetilde\tau_R$ is a smooth function of $(x_2,\xi_1,\xi_2) $,  the energy is propagated along the bicharacteristics, i.e. along the integral curves of the following system of ODEs:
$$ 
  \left\{ \begin{array}{llll}
\displaystyle \dot x^t = \displaystyle  {\nabla_{\xi } \widetilde\tau_R }  (\xi_1^t,x_2^t,\xi_2^t) , & \qquad x^0   =  (x_1^0,x_2^0)   \\
\displaystyle\dot \xi^t= - \displaystyle   {\nabla_x \widetilde\tau_R }  
  (\xi_1^t,x_2^t,\xi_2^t) ,& \qquad \xi^0 =  (\xi_1^0,\xi_2^0) . 
\end{array} \right. $$
Since the condition (\ref{cond1}) avoids the set $ \{ \xi_1 = 0 \} $, we can suppose that $ \xi_1 ^0
\not = 0 $. Moreover since~$ \widetilde\tau_R $ does not depend on $ x_1 $, we find 
$\xi_1^t \equiv \xi_1^0  $. The ODE to be studied is therefore
\begin{equation}\label{bicar}
\left\{ \begin{array}{ll}
\dot x_1^t \! \! \! & \displaystyle
= \, \bar u_1(x_2^t) + { b'(x^t_2) \, (-\xi_1^2+{\xi_2^t}^2+b^2 (x_2^t))\over (\xi_1^2 +{\xi_2^t}^2+b^2 
(x_2^t))^2}   \\
\dot x_2^t \! \! \! & \displaystyle = \, {-2  b'(x_2^t) \, \xi_1{\xi_2^t} \over (\xi_1^2+{\xi_2^t}^2+b^2 (x_2^t))^2} 
  \\
\dot \xi_2^t \! \! \! & \displaystyle = - \bar u_1'(x_2^t)\xi_1 + {2b  ({x_2^t}) b'(x_2^t)^2 \, \xi_1 \over 
(\xi_1^2+{\xi_2^t}^2+b^2 (x_2^t))^2} - {b''(x_2^t) \, \xi_1   \over \xi_1^2+{\xi_2^t}^2+b^2 (x_2^t)} 
  \, \cdotp
\end{array} \right.
\end{equation}
Due to the assumptions on~$\bar u_1$ and on~$b$,  the map $ (x ,\xi ) \mapsto (\nabla_\xi 
\widetilde \tau_R , -\nabla_x  \widetilde \tau_R) (x, \xi) $ is bounded, so   the integral curves  are globally defined in time. The strategy to study their qualitative behaviours 
is to first (in Section \ref{x2xi2}) consider the motion in the reduced phase space $(x_2,\xi_2) \in \R^2 $
and then (in Section \ref{propax1}) to study the motion in the~ $x_1$ direction.


\subsection{Trajectories in the reduced $  (x_2,\xi_2) $ phase space}\label{x2xi2}
In this section  we study the trajectories in the reduced~$ (x_2,\xi_2)$ phase space.  We shall denote~$\xi_1:=\xi_1 ^0$.
\subsubsection{Energy surfaces} 
 Since the Hamiltonian $\widetilde\tau_R$ and $\xi_1$ are  conserved along any trajectory, trajectories are submanifolds of 
 $$ \mathcal E_{\tau,\xi_1} := \bigl \lbrace (x_2,\xi_2) \in \R^2 \, ; \widetilde\tau_R (\xi_1,x_2,\xi_2) = \tau \bigr 
\rbrace . $$
In the following we shall note for any energy~$\tau$ and any $\xi_1\in \R^*$
$$
V_{\tau,\xi_1}(x_2) :=  \frac{ b'(x_2)\xi_1 }{\tau -  \bar u_1 (x_2)\xi_1}  -  \xi_1^2 - b^2 (x_2),
$$
so that if~${\mathcal D} := \displaystyle \bigl\{     x_2 \big / V_{\tau,\xi_1}(x_2) \geq 0\bigr\}
$, then~$\displaystyle \mathcal E_{\tau,\xi_1}   = \Bigl\{
\bigl( x_2,\pm \sqrt{V_{\tau,\xi_1}(x_2)}  \, \bigr) , \: x_2 \in {\mathcal D}
\Bigr\}$. Note that~$V_{\tau,\xi_1}(x_2^t) $ becomes singular if~$x_2^t$ reaches a point~$x_2$ such that~$\tau =  \bar u_1 (x_2)\xi_1$.
\begin{Prop}\label{descriptionenergysurfaces}
The projection of $ \mathcal E_{\tau,\xi_1} $ on the $x_2$-axis is bounded.
\end{Prop}
\begin{proof}
We recall that on~$\mathcal E_{\tau,\xi_1}$ we have
$$
 \frac{ b'(x_2)\xi_1 }{\xi_1^2 + \xi_2^2+ b^2 (x_2)}  + \bar u_1 (x_2)\xi_1 = \tau.
$$
Suppose the trajectory in~$x_2$ is not bounded, then in particular it escapes the support of~$ \bar u_1$. In the case when~$\tau \neq 0$, letting  $ \vert x_2 \vert $ go to infinity yields a contradiction due to the assumptions on~$b$. In the case~$\tau = 0$, for~$x_2$ out of the support of~$ \bar u_1$ we have
$$
 \frac{ b'(x_2)\xi_1 }{\xi_1^2 + \xi_2^2+ b^2 (x_2)} = 0
$$
and the only possibility is for~$x_2$ to be fixed on a zero point of~$b'$ (hence in particular does not go to infinity).
\end{proof}
Proposition~\ref{descriptionenergysurfaces}
 shows that to prove that some trajectories are trapped  in physical space, it suffices to study their behaviour in the~$x_1$ direction. However before doing so, let us prepare that study by  classifying the trajectories in the reduced phase space. 
 Up to a change of parameter, namely expressing time $t$ as a function of $x_2$
 $$ dt =  \pm { b'(x_2)\xi_1\over (\tau-\xi_1 \bar u_1(x_2))^2\sqrt{V_{\tau,\xi_1}(x_2)} } dx_2\,,$$
 which is justified locally, and will give the convenient global behaviour by suitable gluing, we are brought back to the study of the hamiltonian system $\xi_2^2 -V_{\tau,\xi_1}(x_2)$ describing the motion of a particle in the potential $-V_{\tau,\xi_1}$.
 
 For smooth potentials $V$ such that $V(x)\to -\infty $ as $|x| \to \infty$, the possible behaviours of such a system are well-known and the trajectories are usually classified as follows (see for instance~\cite{arnold,AP,BR,hale-kocak,HW}):  periodic orbits,  fixed points, homoclinic and heteroclinic orbits connecting unstable fixed points.
 Here the situation is more complex insofar as $V_{\tau,\xi_1}$ admits singularities. We shall classify the trajectories according to their motion in the~$x_2$ variable.

\subsubsection{Periodic trajectories}\label{defperiodictrajectories}
These  correspond  to the case  when there exists~$[x_{min},x_{max}]$ in~$\R$ with $x_{min}\neq x_{max}$, containing~$x_2^0$ such that 
 \begin{itemize}
 \item $V_{\tau,\xi_1}$ has no singularity and does not vanish on~$]x_{min},x_{max}[$;
 \item $V_{\tau,\xi_1} (x_{min}) =V_{\tau,\xi_1} (x_{max}) = 0 $;
 \item the points~$x_{min}$ and~$x_{max}$ are reached in finite time.
\end{itemize}
The extremal points $x_{min}$ and $x_{max}$ are then turning points, meaning that the motion is periodic. 
The fact that  $x_{min}$ and~$x_{max}$ are reached in finite time is equivalent to 
$V_{\tau,\xi_1}'(x_{min})\neq 0$ and $V_{\tau,\xi_1}'(x_{max})\neq 0$. Indeed if $V_{\tau,\xi_1}'(x_{max})=0$ (resp. $V_{\tau,\xi_1}'(x_{min})=0$), then  $(x_{max},0)$ (resp $(x_{min},0)$) is a fixed point, which contradicts the uniqueness given by the Cauchy-Lipschitz theorem. And conversely, if $V_{\tau,\xi_1}'(x_{max})\neq0$ (resp. $V_{\tau,\xi_1}'(x_{min})\neq 0$), an asymptotic expansion in the vicinity of $x_{max}$ (resp. $x_{min}$) shows that the extremal point is reached in finite time.
\begin{Def}
We will denote by $\PP$ the subset of the phase space~$T^*\R^2$ consisting of initial data corresponding to periodic motions along $x_2$.
\end{Def}
A rather simple continuity argument allows to prove that $\PP$ is an open subset of $\R\times \R^*\times \R^2 $.
Denote indeed by $(\tilde x_1^0,\tilde \xi_1,\tilde x_2^0,\tilde \xi_2^0)$ any point of $\PP$, and by $\tilde x_{min}$ and $\tilde x_{max}$ the extremal points of the corresponding (periodic) trajectory along $x_2$. As $V'_{\tilde \tau,\tilde \xi_1}(\tilde x_{max})\neq 0$ and $V_{\tau,\xi_1}$ is a smooth function of $\tau$ and $\xi_1$ outside from the closed subset of singularity points, the  implicit function theorem  gives the existence of  a neighborhood  of $(\tilde \xi_1,\tilde \tau)$ such that there exists a unique $  x_{max}$ which satisfies $V_{\tau,\xi_1} (x_{max}) =0$. Furthermore $x_{max}$ depends continuously on $\tau$ and $\xi_1$, in particular
$V'_{\tau,\xi_1} (x_{max})\neq 0\,.$
Using the same arguments to build a suitable $x_{min}$, we finally obtain that there exists a neighborhood of $(\tilde x_1^0,\tilde \xi_1,\tilde x_2^0,\tilde \xi_2^0)$ for which the motion along $x_2$ is a non degenerate periodic motion. Moreover, $x_{min}$, $x_{max}$ and also $\xi_{min}$, $\xi_{max}$ and the period $T$ depend continuously on the initial data.

{\it Fixed points} correspond to the degenerate case when $x_{min}=x_{max} =x_2^0$, which implies  that either $\xi_2^0=0$ or $b'(x_2^0)=0$. The latter case is completely characterized by the condition~$b'(x_2^0)=0$, so let us   focus on the case when~$b'(x_2^0)\neq 0$. Fixed points correspond then to local extrema of $V_{\tau,\xi_1}$. They can be either stable or unstable depending on the sign of $V''_{\tau,\xi_1}$.

Stable fixed points are obtained as a limit of periodic orbits when the period $T\to 0$, whereas unstable fixed points are obtained in the limit $T\to \infty$ as explained below.

{\it Stopping trajectories}   belong to the same energy surfaces as unstable fixed points and reach some unstable fixed point in infinite time~:  they 
correspond to the case when  there exists an interval~$[x_{min},x_{max}]$ of~$\R$ containing~$x_2^0$ such that 
 \begin{itemize}
 \item $V_{\tau,\xi_1}$ has no singularity and does not vanish on~$]x_{min},x_{max}[$;
 \item $x_{min}$ and $x_{max}$ are either zeros or singularities of $V_{\tau,\xi_1}$;
  \item as $t\to \infty$, $x_2^t\to x_2^\infty$ such that $V'_{\tau,\xi_1}(x_2^\infty)=0$ or $b'(x_2^\infty)=0$.
\end{itemize}
They are in some sense also a degenerate version of  periodic trajectories since for arbitrarily close initial data, one can obtain periodic orbits. 
\begin{Def}
We will denote by $\delta \PP$ the subset of the phase space~$T^*\R^2$ consisting in initial data corresponding to fixed points and stopping motions along $x_2$.
\end{Def}
Using the characterization of the energy surfaces which carry such pathological motions, we can prove that $\delta \PP$ is a codimension 1 subset of the phase space.
We first consider the energy surfaces containing a fixed point $\tilde x_2$ such that $b'(\tilde x_2)=0$. We then have
$\tau = \xi_1 \bar u_1 (\tilde x_2)$. The corresponding set of initial data
$$\{ ( x_1^0, \xi_1,x_2^0, \xi_2^0) \,/\,  \widetilde \tau_R(  \xi_1,x_2^0, \xi_2^0)= \xi_1 \bar u_1 (\tilde x_2)\}  \hbox{ is of codimension 1.}$$
As we have  assumed  that $b$ has only a finite number of critical points, the union of these sets is still of codimension 1.
We then consider the energy surfaces containing a fixed point $(\tilde x_1, \tilde \xi_1,\tilde x_2,0)$ such that $V_{\tilde \tau,\tilde \xi_1} (\tilde x_2) =V'_{\tilde \tau,\tilde \xi_1} (\tilde x_2)=0$.
We then have
$$ -\bar u_1' (\tilde x_2) +{2b(\tilde x_2)b'(\tilde x_2)\over (\tilde \xi_1^2+b^2(\tilde x_2))^2} -{b''(\tilde x_2) \over \tilde \xi_1^2+b^2(\tilde x_2)}=0.$$
For each $\tilde x_2$, there are at most two values of $\tilde \xi_1^2$ such that the previous quantity vanishes. We therefore deduce that
 $$\displaystyle \Bigl\{ ( x_1^0, \xi_1,x_2^0, \xi_2^0) \,/\, \widetilde \tau_R(  \xi_1,x_2^0, \xi_2^0)= {b'(\tilde x_2) \xi_1\over \xi_1^2+b^2
(\tilde x_2)} +\bar u_1(\tilde x_2) \, \xi_1 \quad \hbox{ and } $$
$$
\qquad \qquad\qquad \qquad\qquad \qquad -\bar u_1' (\tilde x_2) +{2b(\tilde x_2)b'(\tilde x_2)\over (\tilde \xi_1^2+b^2(\tilde x_2))^2} -{b''(\tilde x_2) \over \tilde \xi_1^2+b^2(\tilde x_2)}=0 \Bigr\} 
 \hbox{ is of codimension 1.}$$
 It consists indeed of at most eight manifolds, each one of them   parametrized by  the real parameter $\tilde x_2$.
In the sequel, we shall avoid these pathological motions assuming that the initial data is microlocalized outside   $\delta \PP$.
\subsubsection{Asymptotic trajectories}\label{defasymptotictrajectories}
These correspond to the case when there exists an   interval~$[x_{min}, x_{max}]$  of~$\R$ containing~$x_2^0$ such that
\begin{itemize}
\item $V_{\tau,\xi_1}$ has no singularity and does not vanish on $]x_{min}, x_{max}[$;
\item $x_{min}$ and $x_{max}$ are either zeros or singularities of $V_{\tau,\xi_1}$;
\item as $t\to \infty$, $x_2^t\to x_2^\infty$ where $x_2^\infty\in \{x_{min}, x_{max}\}$ is a pole of multiplicity 1 of $V_{\tau,\xi_1}$. For the sake of simplicity, we further impose that 
$b'(x_2^\infty)\neq 0$.
\end{itemize}
This situation therefore corresponds to a motion which is   not periodic.

Depending on the sign of $\xi_1$, an asymptotic trajectory will either encounter a turning point and  then converge asymptotically to the singular point, or converge monotonically to the limiting point.
As $x_2^\infty$ is such that $\bar u_1(x_2^\infty) = \tau/\xi_1$, one has
\begin{itemize}
\item either $\tau \neq 0$ and $x_2^\infty$ belongs to the support of $\bar u_1$,
\item or $\tau =0$ and $\bar u_1(x_2^0) =-b'(x_2^0) /(\xi_1^2 +\xi_2^2 +b^2(x_2^0)) \neq 0$, meaning that  $x_2^0$ belongs to the support of $\bar u_1$. Therefore, 
$x_2^\infty $ is either $\min \{y>x_2^0 \,/\, \bar u_1 (y)=0\} $ or $\max \{y<x_2^0 \,/\, \bar u_1 (y)=0\} $, in particular $x_2^\infty$ belongs to the support of $\bar u_1$.
\end{itemize}
\begin{Def}
We will denote by $\AA$ the subset of the phase space consisting of initial data corresponding to asymptotic motions along $x_2$.
\end{Def}
The same kind of arguments as in the previous paragraph allow to prove that $\AA$ is an open subset of the phase space $\R\times \R^*\times \R^2$. Consider indeed some $(\tilde x_1^0,\tilde \xi_1,\tilde x_2^0,\tilde \xi_2^0)\in \AA$, and the corresponding asymptotic point $\tilde x_2^\infty$. As $x_2^\infty$ is a pole of multiplicity 1 of $V_{\tilde \tau,\tilde \xi_1}$, the implicit function theorem    shows that $V_{\tau,\xi_1}^{-1}$ admits locally a unique zero, which depends continuously on $\tau$ and $\xi_1$. Using further the continuity of the possible turning point, we get that $\AA$ contains a neighborhood of $(\tilde x_1^0,\tilde \xi_1,\tilde x_2^0,\tilde \xi_2^0)$.
Moreover, we can obtain bounds on the expansion (with respect to time) of any compact subset of $\AA$. Here we will focus on the growth of $\xi_2^t$, and proves that it depends continuously on the initial data in $\AA$.
Without loss of generality, we can consider 
the case when the asymptotic point is~$x_{max}$.
 Then we recall that
 $$
  \lim_{t \rightarrow + \infty} \xi_2^t =   \infty \, 
  \quad \mbox {and} \quad \lim_{t \rightarrow + \infty} x_2^t = x_2^\infty,   \quad \mbox {with} \quad  \xi_1 \bar u_1 
( x_2^\infty)= \tau  .
  $$ 
 As $ x $ tends to $ x_2^\infty$, we have (recalling that~$b'(x_2^\infty) \neq 0$)
  $$
 V_\tau (x) \, \sim \, 
-{b'(x_2^\infty) \over \bar u_1 '(x_2^\infty)}(x-x_2^\infty)^{-1}\, . 
$$
This implies that 
$$| \xi_2^t  |^2  \sim -{b'(x_2^\infty) \over \bar u_1 '(x_2^\infty)}(x-x_2^\infty)^{-1}\quad \mbox{and} \quad \displaystyle \dot x^t_2\sim -2b'(x_2^\infty) \xi_1 {\xi_2^t \over |\xi_2^t|^4} \sim 2{|\bar u_1'(x_2^\infty) |^{3/2} |\xi_1| \over |b'(x_2^\infty|^{1/2}}(x_2^\infty-x)^{3/2}\, . $$
By integration, we get
$$
\displaystyle x^t_2   \sim  x_2^\infty
 + C_1 \ t^{-2}  , \qquad \xi^t_2   \sim 
C_2  t   \,
$$
where $C_1$ and $C_2$ depend continuously on $x_2^\infty$, and consequently on the initial data.

{\it Singular  trajectories} are a degenerate version of the asymptotic trajectories above: they correspond to the case when  there exists an interval~$[x_{min},x_{max}]$ of~$\R$ containing~$x_2^0$ such that 
\begin{itemize}
\item $V_{\tau,\xi_1}$ has no singularity and does not vanish on $]x_{min}, x_{max}[$;
\item $x_{min}$ and $x_{max}$ are either zeros or singularities of $V_\tau$;
\item as $t\to \infty$, $x_2^t\to x_2^\infty$ where $x_2^\infty\in \{x_{min}, x_{max}\}$ is either a singularity of order greater than 1 of $V_{\tau,\xi_1}$, or a singularity which is also a zero of $b'$.
\end{itemize}

As previously, it is easy to check that $x_2^\infty$ is necessarily in the support of $\bar u_1$. Furthermore, we have $\bar u'_1(x_2^\infty)=0.$ Indeed, in the vicinity of $x_2^\infty$, we have
$$V_{\tau,\xi_1} (y) \sim {b'(x_2^\infty)\xi_1\over \bar u(x_2^\infty) - \bar  u(y)}$$
which has a singularity of order greater than 1, or a singularity which is also a zero of $b'$ if and only if $\bar u'(x_2^\infty)=0.$
\begin{Def}
We will denote by $\delta \AA$ the subset of the phase space consisting of initial data corresponding to singular motions along $x_2$.
\end{Def}
The previous condition $\bar u'_1(x_2^\infty)=0$ shows that $\delta \AA$ is included in the union of energy surfaces
$$ \mathcal E_{\tau,\xi_1} \hbox{ with } \tau = \bar u_1(y) \hbox{ for some } y \in supp (\bar u_1) \hbox{ such that }\bar u'_1(y) =0\,.$$
As $\bar u_1'$ has only a finite number of zeros in the support of $\bar u_1$, this implies that $\delta \AA$ is a codimension 1 subset of the phase space.

Gathering all the previous results together, we obtain the following
\begin{Prop}\label{classification}
The phase space  $\R\times \R^* \times \R^2$ admits the following partition
$$ \R\times \R^* \times \R^2 = \PP\cup \AA \cup \Sigma$$
where $\PP$ and $\AA$ are  the open sets of initial data giving rise respectively to periodic motions and asymptotic motions along $x_2$, and $\Sigma : =\delta \PP \cup \delta \AA$ is the codimension 1 set of initial data giving rise to pathological motions along $x_2$.

For any compact set $K\subset \AA \cup \PP$, and for any time $T>0$, we further have a uniform bound on the image of $K$ by the flow up to time $T$.
\end{Prop}
\begin{proof}
The first statement  just tells us that all trajectories belong to one of the four categories described above. Indeed, for any initial data $(x_1^0,\xi_1,x_2^0,\xi_2^0)$, one has
$$V_{\tau,\xi_1}(x_2^0) =(\xi_2^0)^2\geq0,\hbox{ and } V_{\tau,\xi_1} (x_2) \to -\infty \hbox{ as }|x_2|\to \infty,$$
so that there exists  an   interval~$[x_{min}, x_{max}]$  of~$\R$ containing~$x_2^0$ such that
 $V_{\tau,\xi_1}$ has no singularity and does not vanish on $]x_{min}, x_{max}[$,
and  $x_{min}$ and $x_{max}$ are either zeros or singularities of~$V_{\tau,\xi_1}$.
The second statement is then a simple corollary of the continuity results established on $\PP$ and~$\AA$.
 \end{proof}
\subsection{Analysis of the trajectories in the~$x_1$ direction: trapping phenomenon}\label{propax1} $ \, $
Proposition~\ref{descriptionenergysurfaces} states that the trajectories are always bounded in the~$x_2$ variable, so it remains to study the~$x_1$ variable. Our aim is to find a set~$\Lambda  \subset \R \times 
\R^* \times \R^2 $ such that any initial data~$ (x_1^0, \xi_1, x_2^0, \xi_2^0 )$ in~$\Lambda$ gives rise to a trapped 
 trajectory, meaning that
 $$\int_0^t \dot x_1^s ds \hbox{ is uniformly bounded for } t\in \R^+\,.$$
\subsubsection{The criterion of capture}
Let us prove the following result.
\begin{Prop}\label{criterion-capture} 
A necessary and sufficient condition for a trajectory with initial data in~$\AA\cup \PP$ to be trapped is
$$\lim_{t\to T} \frac1t\int_0^t \dot x_1^s ds =0,$$
where $T$ denotes the (finite) period of the motion along $x_2$ in the periodic case, and  $T=+\infty$ in the asymptotic case.
 \end{Prop}
\begin{proof} 
We will study separately the different situations described in the previous section, namely the case of periodic and asymptotic trajectories in $(x_2,\xi_2)$.

$\bullet$ {\it In the case of a  periodic motion} in~$(x_2,\xi_2)$  of period~$ T>0 $,    the function $ \dot x_1^t $ is also periodic, with the same period.
Writing
$$
x_1^t = x_1^0 +\int _0^t \left(  \dot x_1^{s} - {1\over T }\int _0^T  \dot x_1^{s'} \: ds' \right) \: ds+ 
{t\over T }\int _0^T \dot x_1^s \: ds
$$
we see that depending on 
the average of $ \dot x_1^{t }$ over~$ [0,T] $, $x_1^t$ is either a periodic function, or the sum of a periodic function
and a linear function. 
It follows that trapped trajectories are characterized by the   
criterion $ \int _0^T \dot x_1^{t } dt  = 0 $.
Note that, depending on the period $T$, the trajectory can explore  a domain in $x_1$ the size of which may be very large.
Nevertheless the continuity statement in Proposition \ref{classification} shows that we have a uniform bound on this size on any compact subset  of $\PP$.

$\bullet$ For asymptotic motions, we  need to check that 
$$\dot x_1^t - \bar u_1(x_2^\infty)  t \hbox{ is integrable at infinity.}$$
 We have indeed
$$ \dot x_1^t  =\bar u_1( x_2^t ) +{ b'( x_2^t) (-\xi_1^2+{\xi_2^t}^2+b^2(  x_2^t))\over (\xi_1^2+{\xi_2^t}^2+b^2
(x_2^t))^2}, $$
which, together with the asymptotic expansions of $x_2^t$ and $\xi_2^t$ obtained in the previous section, implies that
$$ \dot x_1^t = \bar u_1 (x_2^\infty) +O(t^{-2}). $$
It is  then clear that the trajectory is trapped if and only if 
$$\bar u_1(x_2^\infty)=\lim_{t\to \infty} \frac1t\int_0^t \dot x_1^s ds =0\,.$$
Using again the continuity statement in Proposition \ref{classification},  we also get a uniform bound on the size of the time evolution of any compact subset of $\AA$. 
\end{proof}
\begin{Rem}
Note that, in the case of a singular asymptotic motion along $x_2$, the criterion of capture is equivalent to~$\displaystyle \bar u(x_2^\infty) =\tau =0.$
In the case of periodic motions along $x_2$, this criterion - even more complicated - can also be expressed in terms of the initial data.
The map $ x_2^t$ is indeed a smooth bijection from  a time interval $]t_0,t_0+\frac T2[$ to~$ ]x_{min} , x_{max}[ $, so 
$$
 \int_0^T \dot x_1^t \:  dt \, = \, 2 \ \int_{t_0}^{t_0 + T/2}  \dot x_1^t \: dt \, = \, 2 \
\int_{x_{min}}^{ x_{max}} \,  \dot {x_1} _{\mid \, x_2^t = x_2}^t \  (\dot x_2^t )^{-1}_{
\mid \, x_2^t = x_2} \ dx_2 \, . $$
The trajectory is therefore trapped  if and only if
\begin{equation}
\label{critper}
  \int _{x_{min}}^{x_{max}} \Bigl \lbrack \xi^2_1 - \frac{(\tau / \xi_1) \, b'(y)}{2 \bigl( (\tau / \xi_1) - 
\bar u_1 (y) \bigr)^2} \Bigr \rbrack \, V_{\tau,\xi_1} (y)^{-1/2} \, dy = 0  
\end{equation}
 (notice that since we are only looking for a criterion for the function to vanish, we can replace~$\xi_2$ by~$\sqrt{V_{\tau,\xi_1}} $ without discussing the sign). Such a formula would  be useful  to investigate numerically the initial data giving rise to trapped trajectories. We shall also use it when investigating in more detail the case of the betaplane approximation, in Paragraph~\ref{betaplan}.
\end{Rem}
\subsubsection{Exhibiting a subset of $\Lambda$ of codimension 1}
Let us prove the following proposition.
\begin{Prop}\label{lambda-sing}   Let $ \bar u_1 \in C^\infty_c (\R) $ be a function which is not identically positive, with some zero of finite multiplicity. 
Then  the set $ \Lambda_{sing} $ consisting of data in~$\Lambda $ giving rise to singular and trapped trajectories, is nonempty. It contains a submanifold of $ \R \times \R^* \times \R^2 $ which is of 
codimension~$1 $.
\end{Prop}
\begin{proof} Without loss of generality, we can assume that there exist~$ y_1 < y_2 $  where~$\bar u_1$ vanishes, with~$\bar u'_1( y_2)>0$ and~$\bar u_1 (y) < 0$ on~$ ]y_1,y_2[$.
As we want to study trapped asymptotic trajectories, we will restrict our attention to the case when $\tau =0$,
which is a necessary condition for asymptotic  trajectories to be trapped.
Extremal points of the trajectories are then defined in terms of the function
$$V_{0,\xi_1} (y )= \, - \, \frac{b'(y)}{\bar u_1 (y)} \, - \, b(y)^2 -\xi_1^2\,.$$
More precisely, if we introduce the auxiliary function
$$ \displaystyle \varrho (y) \, := \, - \, \frac{b'(y)}{\bar u_1 (y)} \, - \, b(y)^2 \, , \qquad 
y \in \, ]y_1,y_2[  \, , $$
we obtain  turning (or stopping) points  $y_s$ if $\rho(y_s) =\xi_1^2$, and  singular points $y_s$ if $\displaystyle \lim _{y\to y_s}\rho(y)=+\infty$ or equivalently $\bar u_1(y_s)=0$.
By definition of $y_1$ and $y_2$, one has
$$ \displaystyle \lim_{y \rightarrow y_1 +}  \varrho (y) \, = \, + \infty\quad  \mbox{and } \quad \lim_{y 
\rightarrow y_2-} \quad \varrho (y) \, = \, + \infty.
$$
Let us define~$\displaystyle   N \, := \, \max \ \bigl( \, 0 \, ; 
\inf_{y \in \, ]y_1,y_2[} \ \varrho (y) \, \bigr) \in \R_+ \, . $ 
For $ \xi_1 $ such that $ \xi_1^2 \geq N $, we then define
$$ h(\xi_1) \, := \, \sup \ \bigl \lbrace \, y \in \, ]-\infty , y_2[ \, ; \ \varrho (y) \leq \xi_1^2 \, \bigr \rbrace 
\in \, ]y_1,y_2[ \, . $$
We therefore have that
$$ \forall y \in ]h(\xi_1),y_2[,\quad y \hbox{ is neither a turning point nor a singular point.}$$
As $h$ is a decreasing function on~$]-\infty,  -\sqrt N ] $, all $\xi_1\in]-\infty,  -\sqrt N ] $  except a countable number are continuity points.
Choose then some $\tilde \xi_1$ to be a continuity point of $h$ and $\tilde x_2^0\in ]h(\tilde \xi_1), y_2[$.
By continuity of $h$, there exists a neighborhood $\tilde V$ of $(\tilde x_2^0,\tilde \xi_1)$ such that 
$$ \forall (x_2^0,\xi_1)\in \tilde V, x_2^0-h(\xi_1)>0\,.$$
The set
$\displaystyle \lbrace \bigl ( x^0_1, \xi_1, x^0_2,  {(\varrho (x^0_2) - 
\xi_1^2)^\frac12} \, \bigr) \, ; \ (x^0_1, \xi_1, x^0_2) \in \R \times \tilde V \bigr \rbrace  $
is   a 
submanifold of $ \R \times \R^* \times \R ^2 $ having codimension~1. 
Furthermore, for any initial data in this set, we have $\dot x_{2|t=0} >0$ and a simple connexity argument shows that $x_2^t$ is an increasing function of time. In particular $x_2^t \to y_2$ as $t\to \infty$.
This proves Proposition~\ref{lambda-sing} 
\end{proof}


\subsection{Some examples in the betaplane approximation}\label{betaplan} $ \, $
We are concerned here with 
 the betaplane approximation, that is when $b(x_2) \equiv \beta 
x_2 $. One has therefore 
$$
 \widetilde \tau_R (\xi_1,x_2,\xi_2) = \bar u_1 (x_2) \xi_1 +\frac{\beta \xi_1}{\xi_1^2 + \xi_2^2 +\beta^2 x_2^2}
$$
and
with the notation of Paragraph~\ref{x2xi2},
$$
V_{\tau,\xi_1} (x_2) =  \frac{\beta \xi_1}{\tau} - \xi_1^2 - \beta^2 x_2^2.
$$
$\bullet$ In the absence of convection, one can characterize exactly the set $\Lambda$ of initial data giving rise to trapped trajectories. One can notice that  for~$ \tau $ such that $ \tau \, \xi_1^{-1} \in \, ]0, \beta \, \xi_1^{-2} ] $, the energy surface 
$ \mathcal E_{\tau,\xi_1} $ is simply the ellipse 
$$ \bigl \lbrace \, (x_2 , \xi_2) \, / \, \xi_2^2 + \beta^2 \, x_2^2 =  \beta \, \xi_1 \, \tau^{-1} - 
\xi_1^2 \, \bigr \rbrace \, . $$ 
Let us go through the previous analysis and study the trajectories  in this situation. 
One notices that {\it fixed points} correspond to~$x_2(t) = \xi_2(t) = 0$, with~$x_1(t) = x_1^0$. 
There are no  {\it asymptotic trajectories} (singular trajectories would correspond to~$\tau = 0$, which is not possible here since~$\xi_1 \neq 0$). 
Finally let us consider  {\it periodic trajectories}. In order to get trapping one must check that  if the energy level~$\tau$ and the frequency~$\xi_1$ are fixed, we have
$$
 \int_{x_{min}}^{x_{max}}{\beta \left( {\displaystyle \beta \, \xi_1\over \displaystyle \tau} -2 \, \xi_1^2\right) \over  \sqrt{ 
{\displaystyle\beta\xi_1\over \displaystyle\tau}-\xi_1^2-\beta^2 \, x_2^2}} \ dx_2 = 0  .
$$
Let us choose for instance any~$(x_1,\xi_1) \in \R \times\R^*$ and define the energy level
$ \tau = \beta / (2 \xi_1) $.
Then  the integral is identically equal to zero hence the corresponding trajectory is trapped. It corresponds 
to~$x_1^t = x_1^0$, and to~$x_2^t,\xi_2^t$ satisfying for all times
$ {\xi_2^t}^2 + \beta^2 {x_2^t}^2 = \xi_1^2. $

Such an explicit characterization cannot be obtained if there is some convection, in particular in physically relevant situations.
Nevertheless, we are able to prove that, under suitable assumptions, $\Lambda$ is not empty, and more precisely that  it contains both singular trapped trajectories and periodic 
trapped trajectories.

$\bullet$ Proposition \ref{lambda-sing} shows that under the assumption that
 $\bar u_1 \in C^\infty_c (\R)$  is a function which is not identically positive, with some zero of finite multiplicity,  
$\Lambda$ contains a subset of codimension~1 of initial data giving rise to singular trapped trajectories.

$\bullet$  Finally let us construct periodic trajectories. We suppose to simplify the notation that~$\beta = 1$, and that~$\bar u_1$ has a local maximum at zero, with, say 
$$
0 < \bar u_1 (0) < 2  / 3    \qquad  \mbox{and} \qquad \bar u''_1 ( 0) < - 6  .
$$
Then we shall prove that the set $ \Lambda_{per} $ is nonempty and contains  a submanifold of $ \R \times 
\R^* \times \R^2 $ which is of codimension~$1 $. 
Let~$\tau$ be an energy level. We notice that for~$\eta$ small enough and~$y \in ]-\eta,\eta[$
$$ 
\bar u_1 (y ) <  \bar u_1 (0) <\frac{\tau  }{ \xi_1}  
$$
so there is no singular point for~$V_{\tau,\xi_1}$ in the interval $ ] -\eta, \eta[ $.  
Now  let us show 
that there are  two turning points inside $ ] -\eta, \eta[ $ for~$\eta$ small enough; this   will imply that there is a periodic trajectory of energy~$ \tau  $.
Define the function
$$
H_{\xi_1} (y ) := \frac{\tau  }{\xi_1} - \frac{1}{ \xi_1^2 +y^2} - \bar u_1 (y ) .
$$
Then $ H'_{\xi_1} ( 0) = 0 $  and  an easy computation shows that if~$\eta$ is chosen small enough, then  for all~$(y,\xi_1) \in \, ] -\eta, \eta[ \times [1,+\infty[$,
$$
H''_{\xi_1} (y) = \frac{2 \xi_1^2 - 8y^2 }{ (\xi_1^2 +y^2)^3} - \bar u_1'' (y )\geq 1 / 2  > 0 .  
$$
Now let us choose~$ \eta $ such that 
$\displaystyle \sup_{]-\eta,\eta[} \bar u''_1 (y ) < - 3  ,$
and let us define  the number~$\displaystyle \delta := 1 - \frac{\eta^2}{8 } $. 

We choose from now on~$\xi_1>0$ and  
$\displaystyle
  \tau \equiv \tau 
(\xi_1) := \bar u_1 (0) \xi_1 + \frac{\delta }{ \xi_1}  
$
so that~$\displaystyle H_{\xi_1} (0) =  -   \frac{1-\delta}{ \xi_1^2 } < 0 $. The function $ H_{\xi_1} $ is   decreasing on~$ ] -\eta, 0[ $  and  increasing on the interval $ ] 0,  \eta[ $; moreover
$$ H_{\xi_1} (  \pm \eta) \geq H_{\xi_1} ( 0) + \frac{\eta^2}{4 } = - (1-\delta) + \frac{\eta^2}{4 } 
\geq  \frac{\eta^2}{8  } > 0 . $$
It follows that, for all $ \xi_1 \in [1,+\infty[ $, there are two points~$x_{min}(\xi_1)$ (in~$] -\eta, 0[$) and~$x_{max}(\xi_1)$  (in~$] 0,  \eta[$) satisfying the requirements of periodic trajectories, in the sense of Section~\ref{defperiodictrajectories}.

According to the criterion (\ref{critper}), we define the following function on~$ [1,+\infty[$: 
$$ 
G( \xi_1 ) := \int _{x_{min}(\xi_1)}^{x_{max}(\xi_1)} \Bigl \lbrack \xi^2_1 - \frac{(\tau / \xi_1
)  }{2 \bigl( (\tau / \xi_1) - \bar u_1 (y) \bigr)^2} \Bigr \rbrack \, V_{\tau,\xi_1} (y)^{-1/2} \, dy ,
$$
and let us prove it vanishes.
We notice that recalling that~$ \bar u_1 ( 0)  < 2/3$,
\begin{eqnarray*}
G(1)& \geq& \int _{x_{min} (1)}^{x_{max} (1)} \Bigl \lbrack 1 - \Bigl( \bar u_1 ( 0) 
+ \delta \Bigr) \frac{1}{2 \delta^2} \Bigr \rbrack \, V_{\tau,\xi_1}  (y)^{-1/2} \, dy\\
& \geq &
\frac{6 \delta^2 -2 -3 \delta}{6 \delta^2} \int _{x_{min} (1)}^{x_{max} (1)}V_{\tau,\xi_1} 
(y)^{-1/2} \, dy > 0 
\end {eqnarray*}
since~$\eta\in ]0,1[$. 
On the other hand, when $ \xi_1 $ goes to $ + \infty $, we find that $x_{min} (\xi_1) \rightarrow  0 $
and $x_{max}(\xi_1) \rightarrow  0 $ so that recalling that~$\tau (\xi_1) = \bar u_1 ( 0) \xi_1 - \delta /\xi_1$,
$$ G( \xi_1 ) \sim - \bar u_1 ( 0) \, \frac{ \xi_1^4}{2  \delta^2} \int _{x_{min}(\xi_1)}^{x_{max}
(\xi_1)} V_{\tau,\xi_1}  (y)^{-1/2} \, dy , $$
so~$\displaystyle  \lim_{\xi_1 \rightarrow + \infty} \ G( \xi_1 ) <0 .  $
By construction, the function $ G $ is smooth so there is some $ \tilde \xi_1$ belonging to~$]1,+\infty[ $ such that 
$ G(\tilde \xi_1) = 0 $. We have generically (in~$\tau$) $ G' (\tilde \xi_1) \not = 0 $. This implies that~$ \Lambda_{per} \not = \emptyset $, where~$ \Lambda_{per}$ is the subset of~$\Lambda $ giving rise to periodic trajectories,  and that~$ \Lambda_{per}  $ contains a submanifold of 
$ \R \times \R^* \times \R^2 $ which is of codimension~1.

\section{Study of the Poincar\'e waves}\label{sectionpoincare}
\subsection{The strategy}
In this section we want to prove the dispersion property of the Poincar\'e polarisation, namely result~(2) of Theorem~\ref{mainresult}.
We recall that the principal symbols~$\tau_\pm$ are of order one, so
one needs to study diffractive-type propagation, on a time  scale of the order~$  1 / \varepsilon^2$. 
Computing the bicharacterstics, as in the previous section, is therefore not enough to understand the classical flow. We shall instead rely on a spectral argument to prove that  Poincar\'e waves do disperse, and escape from any compact set in the physical
space.
The first step consists in taking the Fourier transform in~$x_1$ (recalling that the problem is invariant by translations in the~$x_1$ direction). We also recall that the data is microlocalized on a compact set~${\mathcal C}$ such that~$\xi_1$ is bounded away from zero.
Since  $\tau_\pm=\pm\sqrt{\xi_2^2+b^2(x_2)+\xi_1^2 } $,    functional 
calculus implies that one can find classical  pseudo-differential operators $H_{ \pm}(\xi_1)$ of principal symbols 
$\xi_2^2+b^2(x_2)$ such that $T_\pm=\pm\sqrt{H_{\pm}(\xi_1)+\xi_1^2 }$.  
Let us now  call $\lambda_\pm^k(\xi_1)$ and $\varphi_\pm^k(\xi_1;x_2)$ the eigenvalues and eigenfunctions of $  H_{\pm}(\xi_1)$. 
The following proposition will be proved in the next paragraph.
\begin{Prop}\label{cinfty}
Let~$\phi$ be an eigenfunction
 of $  H_{\pm}(\xi_1)$, microlocalized on an energy surface which interstects~${\mathcal C}$.
Then~$\phi$ and its associate eigenvalue~$\lambda$ are~$C^\infty$ functions of~$\xi_1$. Moreover~$\frac1 \e \partial_{\xi_1} \lambda$ is bounded on compact sets in~$\xi_1$.
\end{Prop} 
\subsection{Proof of Theorem~\ref{mainresult}(2) assuming  Proposition~\ref{cinfty}}
Let us now carry out this program. We consider an initial data denoted~$\varphi^0$, microlocalized in~${\mathcal C}$. One can take the Fourier transform in~$x_1$ which gives
$$
\varphi^0(x) = \frac1{\sqrt{2\pi\e}} \int \hat \varphi^0(\xi_1,x_2) e^{-i\frac {x_1 \xi_1}{\e}} \: d\xi_1.
$$
Now let us consider a    coherent state (in Fourier variables) at $(q,p)$ (see Appendix~\ref{semiclassic}), that is:
\be
\varphi_{qp}(\xi_1):=  \frac1{ (\pi\eps)^{ \frac{1}4}}  e^{i\frac{\xi_1 q}{\eps} } e^{-\frac{(\xi_1-p)^2}{2\eps}}.
\ee
After decomposition onto coherent states we get
$$
\varphi^0(x) = \frac1{\sqrt{2\pi\e} } \frac1{ (\pi\eps)^\frac14}\int \tilde \varphi^0(q,p,x_2) e^{i\frac{\xi_1 (q-x_1) }{\eps} } e^{-\frac{(\xi_1-p)^2}{2\eps}}\: dq dp d\xi_1 
$$
where
$\displaystyle
 \tilde \varphi^0(q,p,x_2)  := \left( \varphi_{qp}  |  \hat \varphi^0(\cdot,x_2) \right)_{L^2}.
 $
We notice that the integral over~$p$ and~$q$ is, modulo~$O(\e^\infty)$, on a compact domain due to the microlocalization assumption on~$\varphi^0$. Finally decomposing onto the eigenfunctions~$\varphi_\pm^k(\xi_1;x_2)$ gives
$$
\varphi^0(x) = \frac1{\sqrt{2\pi\e} } \frac1{ (\pi\eps)^\frac14} \sum_k \int \overline \varphi^0(q,p;k,\xi_1,x_2) \varphi_\pm^k(\xi_1;x_2) e^{i\frac{\xi_1 (q-x_1) }{\eps} } e^{-\frac{(\xi_1-p)^2}{2\eps}}\: dq dp d\xi_1 
$$
where
$
\displaystyle \overline \varphi^0(q,p;k,\xi_1,x_2) := \left(\varphi_\pm^k(\xi_1;\cdot)  |  \tilde \varphi^0(q,p,\cdot)  \right)_{L^2}.
$
Note that the dependence of~$\overline \varphi^0$ on~$\xi_1$ is only through the eigenfunction~$\varphi_\pm^k$, so~$\overline \varphi^0$ depends smoothly on~$\xi_1$, as stated in Proposition~\ref{cinfty}.

The sum over~$k$ contains~$O(\e^{-1})$ terms, due to the fact that~$\lambda_\pm^k(\xi_1)$ remains in a finite interval (this will be made more precise in the next section, see Remark~\ref{supportpoincare}). Now it remains to propagate  at time $t  /\varepsilon^2$ this initial data, 
which  gives rise to the following expression:
$$
\frac1{\sqrt{2\pi\e} } \frac1{ (\pi\eps)^\frac14} \sum_k \int \overline \varphi^0(q,p;k,\xi_1,x_2) \varphi_\pm^k(\xi_1;x_2) e^{i\frac{\xi_1 (q-x_1) }{\eps}} e^{-\frac{(\xi_1-p)^2}{2\eps}} 
e^{\pm i (\lambda_\pm^k(\xi_1) +\xi^2_1)^{\frac 1 2} \frac t{\e^2} }
 \: dq dp d\xi_1 
$$
The stationary phase lemma then  gives that this integral is $O(\varepsilon^\infty)$ except if there exists a stationary point,  given by the conditions:
\[
\xi_1=p\ \mbox{and}\  \varepsilon(x_1-q)\pm\frac{ 2\xi_1 +\partial_{\xi_1}\lambda_\pm^k(\xi_1) }{2\sqrt{\lambda_\pm^k(\xi_1)+\xi^2_1}}t=0 . \]
The second condition gives 
$$
2\sqrt{\lambda_\pm^k(\xi_1)+\xi^2_1}(x_1-q)=\mp \frac1\varepsilon\left( 2p+\partial_{\xi_1}\lambda_\pm^k(\xi_1)\right)t.
$$
Therefore, since $p\neq 0$ and the $\lambda^k_\pm$\ 's are bounded, with~$\partial_{\xi_1}\lambda_\pm^k(\xi_1) = O(\e)$, there is no critical point for $x_1$ in a compact set.
 Proposition~\ref{cinfty} therefore
allows to   apply the stationary phase lemma and to conclude the proof of  result~(2) of Theorem~\ref{mainresult}. 
Notice that the (fixed) losses in~$\e$ (namely the negative powers of~$\e$ appearing in the integrals and the number of~$k$'s in the sum) are compensated by the fact that the result is~$O(\e^\infty)$; it is important at this point that as noticed above,  the function~$\overline \varphi^0$ depends smoothly on~$\xi_1$. 
\subsection{Proof of Proposition~\ref{cinfty}}
The first step of the proof consists in using the theory of normal forms in order to reduce the problem to the study of functions of the harmonic oscillator (as in Bohr-Sommerfeld quantization, paying special attention to the dependence on~$\xi_1$). The second step then consists in checking that the eigenvectors and eigenvalues have the required dependence on~$\xi_1$. 
In the following we shall only deal with~$H_+$ to simplify, and we shall write~$H := H_+$.
The first step relies on the following lemma. 
\begin{Lem}
Suppose that~$\xi_1$ lies in a compact set away from zero. There is an elliptic Fourier Integral Operator~$U : {\mathcal D}(\R) \to {\mathcal D}(\R)$, independent of~$\xi_1$, and a pseudodifferential operator~$V(\xi_1)$ with~$C^\infty$ symbol in~$\xi_1$ such that microlocally in any compact set~$K \subset T^*\R$ one has~$(VU)^* = (VU)^{-1}$ and
$$
VUH (VU)^{-1} = f(-\e^2 \partial_2^2 +x_2^2 ; \xi_1;\e) + O(\e^\infty)
$$
where
$$
 f(I ; \xi_1;\e) = f_0(I) + \sum \e^j f_j (I ; \xi_1;\e)
$$
while~$f_0$ is a smooth bijection on~$\R^+$ and the~$f_j$'s are~$C^\infty$ functions of~$I$ and~$\xi_1$.
\end{Lem}
\begin{proof}
The proof consists in using  techniques linked to the isochore Morse lemma (see~\cite{colinvey}, \cite{san}). Actually
we introduce a canonical change of variables (the corresponding operator being  the  FIO~$U$) allowing to pass from the variables~$(x_2,\xi_2)$ to action-angle variables (see~\cite{arnold} for instance). Let us make this first step more precise : we recall that the action variables are given by
$$
I:= \oint \xi_2 \: dx_2,
$$
where the integral is taken on a constant energy curve~$H = h$. The angle variables are then given by solving~$  \theta = \partial_I S (I, \xi_2)$, where~$K$ is the hamiltonian in the new variables, which only depends on~$I$, and where~$S$ is defined by 
$$
dS_{I = constant} = \xi_2 dx_2.
$$
Note that this is a {\it global} change of variables.
It is now well known (see~\cite{DS} for instance) that such a canonical change of variables is associated with   an FIO~$U$ (independent of~$\xi_1$ since the principal symbol of~$H$ does not depend on~$\xi_1$) such that
\begin{equation}\label{uhu}
UHU^{-1}  = f^0(-\e^2 \partial_2^2 +x_2^2) + \e F_1
\end{equation}
where~$f^0$ is a smooth, global bijection on~$\R^+$.  In this case we can actually write (see~\cite{DS}) the following formula for~$U$:   for any~$L^2$ function~$\varphi$
$$
U \varphi(x_2) =\frac1{(2\pi \e)^\frac32} \int e^{i\frac{S(x_2,\xi_2 )- y \xi_2}{\e}} a_0(x_2,\xi_2,\e) u(y)\: dyd\xi_2
$$
where~$a_0$ is constructed so that~$U$ is unitary (up to~$O(\e)$).
Using those new coordinates and the formula giving~$U$ as well as the formula for the principal symbol of the adjoint (which is here the inverse) given in Appendix~\ref{semiclassic}, it is not difficult to show that~(\ref{uhu}) holds.
Once the function~$f^0$ is obtained we proceed by induction: we first look for a symbol~$q_1$ such that~$Q_1:= {\rm Op}_\e^0(q_1)$ satisfies
$$
e^{iQ_1} UHU^{-1} e^{-iQ_1} =  f^{0,1}(-\e^2 \partial_2^2 +x_2^2) + \e^2 F_2.
$$
In order to compute~$q_1$ we notice that
$$
e^{iQ_1} UHU^{-1} e^{-iQ_1} =  UHU^{-1} + \int_0^t \partial_t H_t \: dt
$$
where
$\displaystyle
 H_t:= e^{itQ_1} UHU^{-1} e^{-itQ_1} .
$
One sees easily that
\begin{eqnarray*}
\partial_t H_t &=&  i [Q_1,H_t]\\
 &=&e^{itQ_1} [Q_1,  f^0(-\e^2 \partial_2^2 +x_2^2) + \e F_1]e^{-itQ_1}
\end{eqnarray*}
so the principal symbol of~$\partial_t H_t$ is therefore the Poisson bracket
$
\displaystyle \e \{q_1,f^0(x_2^2 + \xi_2^2)\}.
$
One then remarks that the equation
$$
\{q,f^0(x_2^2 + \xi_2^2)\} = g
$$
has a solution if and only if~$g$ has zero mean value (in action-angle variables): one has indeed necessarily
$$
\partial_\theta q = \frac{I g}{\partial_I f} \cdotp
$$
Note that if~$f$ and~$g$ are smooth, then so is~$q$ since~$\partial_I f >0$ everywhere. It follows that one can find~$q_1$ simply by solving
$$
\displaystyle   \{q_1,f^0(x_2^2 + \xi_2^2)\} = f_1 - \overline f_1
$$
where~$ \overline f_1$ is the average of~$f_1$
$$
 \overline f_1 (I) = \frac1{2\pi} \int_0^{2\pi} f_1(R_\theta (x_2,\xi_2) ; \xi_1) \: d\theta
$$
where~$R_\theta$ denotes the rotation of angle~$\theta$. This implies that with this choice of~$q_1$ and writing~$Q_1:=  {\rm Op}_\e^0(q_1)$ one has
$$
e^{iQ_1} UHU^{-1} e^{-iQ_1} =   (f^{0 } + \e \overline f_1)(-\e^2 \partial_2^2 +x_2^2) + \e^2 F_2.
$$
One proceeds similary at all orders.
\end{proof}
The next step consists in using that lemma to check that the eigenvalues and eigenfunctions enjoy the expected smoothness properties. Actually this is rather straightforward since if~$(\psi, \lambda)$ satisfy
$$
H\psi = \lambda \psi + O(\e^\infty)
$$
then~$\widetilde \psi := VU \psi$ satisfies
$$
VUH (VU)^{-1} \widetilde \psi = \lambda \widetilde \psi   + O(\e^\infty)
$$
hence
$$
 f (-\e^2 \partial_2^2 +x_2^2)  \widetilde \psi = \lambda \widetilde \psi   + O(\e^\infty) .
$$
This implies that
$$
\lambda = \lambda_n = f\left( {\e} \bigl(n + \frac12 \bigr) ; \xi_1 ; \e\right) \quad \mbox{and} \quad  \widetilde \psi  =  \widetilde \psi_n = \frac1{\e^\frac14} h_n\left( \frac{x_2}{\sqrt \e}\right) 
$$
where~$h_n$ is the $n$-th Hermite function. 
The conclusion follows recalling that~$ f(I ; \xi_1;\e) = f_0(I) + \sum \e^j f_j (I ; \xi_1;\e) $ (where~$f_0$ does not depend on~$\xi_1$) and that each eigenfunction is microlocalized in a compact set.
\begin{Rem}\label{supportpoincare}
Since we have a complete spectral description of~$T_\pm$, with discrete spectrum for each given~$\xi_1$, it is obvious that if the initial data is microlocalized in~$(\xi_1,x_2,\xi_2)$, then the solution to the equation~$\e^2 \partial_t \varphi =i T_{\pm} \varphi$ remains  microlocalized in  the set of energy surfaces containing~$(\xi_1,x_2,\xi_2)$. Notice also that there are~$O(1/\e)$ eigenvalues in a compact energy surface.
\end{Rem}

 \section{Diagonalization}\label{diago}

 In this section we shall prove that the scalar propagators defined in Section~\ref{scalar} correspond indeed to a diagonalization of the original linear system~(\ref{linsyst}).  Let us prove the following proposition. 
 \begin{Prop}
 Consider a compact set~${\mathcal K}   \subset \R^* \times T^*\R  $. 
  With the notation of Definitions~\ref{therossbyproj} and~\ref{thepoincareproj},  the operator~$\Pi:=( \Pi_- \, \Pi_R \,  \Pi_+ )$ maps continuously~$H^\infty_\e$ functions, microlocalized in~$\R \times {\mathcal K}$ and satisfying~(\ref{cond1}), onto~$H^\infty_\e$. Moreover it is left-invertible modulo $\varepsilon^\infty$, and its   left inverse~$Q$ (modulo $\varepsilon^\infty$)   maps continuously those functions onto~$H^\infty_\e$.
 \end{Prop}
 \begin{Rem}
 This proposition, along with Proposition~\ref{classification} and Remark~\ref{supportpoincare}  showing the propagation of the microlocal support (in~$(\xi_1,x_2,\xi_2)$) of Rossby and Poincar\'e modes, prove the first part of Theorem~\ref{mainresult}. 
 \end{Rem}
 \begin{proof} The main step consists in showing that~$\Pi$ does have a left inverse, and in computing its principal symbol.
 The construction of a left inverse can be done symbolically as follows.
We first compute the matrix-principal symbol $\mathcal P_0$ of $\Pi $. One gets 
$$
\mathcal P_0=
\left(\begin{array}{ccc}
\displaystyle \frac{-\xi_1\sqrt{\xi^2 + b^2} - i\xi_2 b }{\xi_1^2 + b^2}
 & \displaystyle-\frac{ib}{\xi_1}
& \displaystyle \frac{ \xi_1\sqrt{\xi^2 + b^2} - i\xi_2 b }{\xi_1^2 + b^2}
\\
1& \displaystyle - {\xi_2 \over \xi_1}
&1\\
\displaystyle  \frac{\xi_1\xi_2 +i b\sqrt{\xi^2 + b^2}  }{\xi_1^2 + b^2}
&1&
\displaystyle  \frac{\xi_1\xi_2 -i b\sqrt{\xi^2 + b^2}  }{\xi_1^2 + b^2}
\end{array}\right).
$$
This shows that~$\Pi$ maps microlocalized functions in~$\R \times {\mathcal K}$ onto~$H^\infty_\e$. 
A simple computation shows that
$$
|\mbox{det} \: \mathcal P_0 |=\frac{ 2 (\xi^2 + b^2)^\frac32}{(\xi_1^2+b^2) |\xi_1|} \geq 2,$$
therefore $\mathcal Q_0 = \mathcal P_0^{-1}$ exists.
Let us call $Q_0$ the matrix obtained by Weyl quantization (term by term) of $\mathcal Q_0 $. By symbolic calculus we have that:
$$
Q_0 \Pi ={\rm Id}+\varepsilon I_1
$$
Let us call $\mathcal Q_1=-\mathcal I_1\mathcal Q_0^{-1}$, where $\mathcal I_1$ is the (matrix) principal symbol of $I_1$. Again by symbolic calculus we have that
$$
\left(Q_0+\varepsilon Q_1\right) \Pi ={\rm Id}+\varepsilon^2I_2
$$
 where $Q_1$ has principal symbol $\mathcal Q_1$.
Defining now $\mathcal Q_2=-\mathcal I_2\mathcal Q_0^{-1}$ we get 
$$
\left(Q_0+\varepsilon Q_1+\varepsilon^2Q_2\right) \Pi ={\rm Id}+\varepsilon^3I_3
$$ 
where $Q_2$ has principal symbol $\mathcal Q_2$,
and so on.
 This allows to invert the matrix~$\Pi$ up to~$O(\e^\infty)$, and the principal symbol of the (approximate) inverse matrix~$Q$ is given by~${\mathcal Q}_0 $: we have
 $$
{\mathcal Q}_0 = \frac{1}{2(\xi^2 + b^2)}
\left(\begin{array}{ccc}
\displaystyle ib\xi_2 - \xi_1 \sqrt{\xi^2 + b^2} &
 \displaystyle \xi_1^2 + b^2 
 & \displaystyle  - ib\sqrt{\xi^2 + b^2} + \xi_1 \xi_2 
\\
\displaystyle 2ib\xi_1  & \displaystyle - 2\xi_1\xi_2  
& \displaystyle 2\xi_1^2  \\
\displaystyle   ib\xi_2 + \xi_1 \sqrt{\xi^2 + b^2}   &
 \displaystyle     \xi_1^2 + b^2   &
 \displaystyle  ib\sqrt{\xi^2 + b^2 } + \xi_1 \xi_2  
 \end{array}\right).
 $$
 Note that each term of the expansion of the symbol of~$Q$ is a polynomial of increasing order in~$\xi$, but that is not a problem due to the microlocalization assumption. The operator~$Q$  therefore
    clearly maps continuously~$H^\infty_\e$ onto itself.
The proposition is proved.
\end{proof}

%
 
\section{Control of the nonlinear terms for a weak coupling}\label{nonlinear}

We are now interested in describing the behaviour of our initial nonlinear system, which includes the effect of the convection by the unknown $ (u_1,u_2)$. We shall prove Theorem~\ref{nonlineartheorem}  in this section.
We recall that the system reads as follows:
\begin{equation}\label{NLS}
\eps^2 \d_t U+ A(x,\eps D_x) U +\eps ^{3+\eta} S_1(U) \eps \d_1 U +\eps^{3+\eta}  S_2(U) \eps \d_2 U=0 , \quad \eta \geq 0
\end{equation}
with
\begin{equation}
 S_1 (U) = \left( \begin{array}{ccc}
u_1 & \frac{1}{2} u_0 & 0 \\  
 \frac{1}{2}u_0 & u_1 & 0 \\
 0 & 0 & u_1  \end{array} \right) , \qquad S_2 (U) = \left( \begin{array}{ccc}
u_2 & 0 & \frac{1}{2}u_0 \\  
0 & u_2 & 0 \\
 \frac{1}{2}u_0 & 0 & u_2 \end{array} \right)  
 \end{equation}
 and~$U = (u_0, u_1, u_2)$.
 The usual theory of symmetric hyperbolic systems   provides the local existence of a solution to
(\ref{NLS}) in $H^s(\R^2)$ for $s>2$, on a time interval depending on~$\e$ a priori.
Because of the semiclassical framework (fixed by the form of the initial data), it is actually natural to rather consider $\eps$-derivatives. Moreover, as the derivative with respect to~$x_2$ does not commute with the singular perturbation, we expect even the semiclassical Sobolev norms to grow like $\exp \left(\frac{Ct}\eps\right)$, and therefore the life span of the solutions to (\ref{NLS}) to be non uniform with respect to $\eps$. That is the reason why the~$W^{s }_\e$  spaces were introduced in~(\ref{defhs1s2}).
\subsection{Propagation of regularity for the linear singular perturbation problem}$ $
Let us first remark that derivatives with respect to $x_1$ do commute with $A(x,\eps D_x)$, so   we can propagate as much regularity in $x_1$ as needed.
Extending a result by Dutrifoy, Majda and Schochet \cite{DU} obtained in the particular case when~$
b(x_2) =\beta x_2$, we will actually prove that there is an operator of principal symbol~$(\xi_2^2+b^2) {\rm Id}$ which ``almost commutes" with $A(x,\eps D_x)$ in the semiclassical regime.
 
$\bullet$ The first step, as in~\cite{DU},  is to perform the following orthogonal change of variable
$$ \tilde U : = \left( {u_0+u_1\over \sqrt{2}}, {u_0-u_1\over \sqrt{2}}, u_2\right)$$
in order to produce the generalized creation and annihilation operators
$$
L_\pm :=\frac1{\sqrt{2}} \bigl( \eps \d_2 \mp b \bigr)\,.
$$
The system (\ref{NLS}) can indeed be rewritten
$$
\eps^2 \d_t \tilde U +\tilde A(x,\eps D_x) \tilde U +\eps ^3 \tilde S_1(\tilde U) \eps \d_1 \tilde U +\eps^3 \tilde S_2(\tilde U) \eps \d_2 \tilde U =0$$
with 
$$
\tilde A(x,\eps D_x):=\begin{pmatrix}
\eps \bar u_1   \eps \partial_1 +\eps \d_1&0& L_+ +\frac{\eps^2}{\sqrt{2} } \bar u_1'\\
0 & \eps \bar u_1   \eps  \partial_1 -\eps \d_1 & L_- -\frac{\eps^2}{\sqrt{2} } \bar u_1'\\
 L_-&L_+&\eps \bar u_1   \eps  \partial_1  
\end{pmatrix}\,,
$$
and
$$\tilde  S_1 (\tilde U) := \left( \begin{array}{ccc}
\displaystyle {3\tilde u_0-\tilde u_1\over 2\sqrt{2}} & 0 & 0 \\  
0 & \displaystyle{\tilde u_0-3\tilde u_1\over 2\sqrt{2}} & 0 \\
 0 & 0 & \displaystyle{\tilde u_0-\tilde u_1\over  \sqrt{2}}  \end{array} \right) , \qquad \tilde S_2 (\tilde U) := \left( \begin{array}{ccc}
\tilde u_2 & 0 & \displaystyle {\tilde u_0+\tilde u_1\over 4} \\  
0 & \tilde u_2 & \displaystyle {\tilde u_0+\tilde u_1\over 4}  \\
\displaystyle{\tilde u_0+\tilde u_1\over 4}  & \displaystyle{\tilde u_0+\tilde u_1\over 4}  & \displaystyle \tilde u_2 \end{array} \right) . 
$$

 $\bullet$ Next, remarking that
$\displaystyle[\eps^2\d_2^2-b^2, \eps\d_2 \pm b] = \pm 2\eps b ' (\eps \d_2 \pm b) \pm \eps^2 b'',$
we introduce the operator
$$D_\eps:=\begin{pmatrix}
\eps^2\d_2^2-b^2+2\eps b '&0&0\\0&\eps^2\d_2^2-b^2-2\eps b'&0\\
0&0&\eps^2\d_2^2-b^2
\end{pmatrix}.$$
We notice that~$D_\e$ is a scalar operator at leading order. Moreover one can
  compute the commutator~$[D_\eps, \tilde A(x,\eps D_x)]$: we find
$$[D_\eps, \tilde A(x,\eps D_x)] =  \begin{pmatrix}
 \displaystyle
[\e^2 \partial_2^2, \e \bar u_1 ]  \eps \d_1  &0& \displaystyle\frac{\e^2}{\sqrt{2}}  \bigl(
[\e^2 \partial_2^2,  \bar u'_1 ]  + 2 \e  b'  \bar u'_1 -   b''
\bigr)\\
0&[\e^2 \partial_2^2, \e \bar u_1 ]  \eps \d_1 &  \displaystyle - \frac{\e^2}{\sqrt{2}}  \bigl(
[\e^2 \partial_2^2,  \bar u'_1 ]  - 2 \e  b'  \bar u'_1 -   b''
\bigr)\\ \displaystyle
- \frac{\eps^2}{\sqrt{2}}b''& \displaystyle \frac{\eps^2}{\sqrt{2}}b''& [\e^2 \partial_2^2, \e \bar u_1 ]  \eps \d_1  
\end{pmatrix}
$$
 
$\bullet$ At this stage, we have proved that
\begin{equation}
\label{commutator}
 [D_\eps, \tilde A(x,\eps D_x)] = O(\eps^2( {\rm{Id}}-\eps^2\d_1^2-D_\eps))
 \end{equation}
 meaning that the commutator~$ [D_\eps, \tilde A(x,\eps D_x)] $ is of order~$O(\e^2)$ with respect to the elliptic operator $  {\rm{Id}}-\eps^2 \partial^2_1 - D_\eps $. That implies that the regularity of the solution to the linear equation 
$$\eps^2 \d_t V+ \tilde A(x,\eps D_x) V=0$$
can be controlled by an application of Gronwall's lemma: one has
$$
\eps^2 \|(Id-\eps^2\d_1^2-D_\eps) V(t)\|_{L^2}^2 \leq \eps^2 \| (Id-\eps^2\d_1^2-D_\eps) V_0\|^2_{L^2}+ C\eps^2\int_0^t \| (Id-\eps^2\d_1^2-D_\eps) V(s)\|^2_{L^2}ds\,,$$
where $C$ depends on the $W^{2,\infty}$ norms of $\bar u_1$ and $b$, so
$$
 \|(Id-\eps^2\d_1^2-D_\eps) V(t)\|_{L^2}^2 \leq C \| (Id-\eps^2\d_1^2-D_\eps) V_0\|^2_{L^2} \: e^{  C t}.
$$
\subsection{Uniform a priori estimates for the nonlinear equation}
Since the extended harmonic oscillator controls two derivatives in $x_2$, we get a control on the Lipschitz norm of $U$ of the type
\begin{equation}
\label{lipschitz-control}
\|\e \partial_j U\|_{L^{\infty} } \leq \frac C\eps( \| D_\eps^2 U \|_{L^2} + \| \eps^4 \d_1^4 U\|_{L^2}+\|U\|_{L^2})\,.
\end{equation}
As $D_\eps$ is a scalar differential operator at leading order in $\eps$, the antisymmetry of the higher order nonlinear term is preserved. More precisely, we have, using the  Leibniz formula, 
$$\begin{aligned}\eps^2 \d_t D_\eps \tilde U &+\tilde A(x,\eps D_x) D_\eps \tilde U  +\eps ^3 \tilde S_1(\tilde U) \eps \d_1 D_\eps \tilde U +\eps^3 \tilde S_2(\tilde U) \eps \d_2 D_\eps \tilde U\\
=-& [D_\eps, \tilde A(x,\eps D_x)] \tilde U - \eps^3 \tilde S_2(\tilde U) [D_\eps ,\eps\d_2]\tilde U - 
\eps^3 [D_\eps ,  \tilde S_j(\tilde U)]  \eps \d_j \tilde U
\end{aligned}  $$
as well as
$$\begin{aligned}\eps^2 \d_t D^2_\eps \tilde U &+\tilde A(x,\eps D_x) D^2_\eps \tilde U  +\eps ^3 \tilde S_1(\tilde U) \eps \d_1 D^2_\eps \tilde U +\eps^3 \tilde S_2(\tilde U) \eps \d_2 D^2_\eps \tilde U = - [D^2_\eps, \tilde A(x,\eps D_x)] \tilde U \\
- & \eps^3 \tilde S_2(\tilde U) [D_\eps ,\eps\d_2]D_\eps\tilde U - 
\eps^3 [D_\eps ,  \tilde S_j(\tilde U)]  \eps \d_j D_\eps\tilde U \\
+ & D_\eps \left( - \eps^3 \tilde S_2(\tilde U) [D_\eps ,\eps\d_2]\tilde U - 
\eps^3 [D_\eps ,  \tilde S_j(\tilde U)]  \eps \d_j \tilde U\right).
\end{aligned}  $$
and in the same way, for~$1 \leq \ell \leq 4$,
$$\begin{aligned}\eps^2 \d_t ( \eps \d_1) ^ \ell \tilde U &+\tilde A(x,\eps D_x)  ( \eps \d_1) ^ \ell \tilde U  +\eps ^3 \tilde S_1(\tilde U) ( \eps \d_1) ^{\ell+1}\tilde U +\eps^3 \tilde S_2(\tilde U)  \eps  \d_2   ( \eps \d_1) ^ \ell \tilde U\\
&=  - \eps^4\sum _{k=1}^ \ell  C^ \ell_4  (\eps\d_1)^ \ell \tilde S_j(\tilde U)\eps\d_{j}(\eps\d_1)^{\ell-k} \tilde U\,.
\end{aligned}  $$
In all  cases, the terms of higher order disappear by integration in~$x$ and the other terms are controlled
 with the following trilinear estimate (writing generically~$ \tilde Q(\tilde U)$ for all the nonlinearities): for all~$0 \leq k \leq 2$ and all~$0 \leq \ell \leq 4$  
\begin{equation}
\label{trilinear}
\begin{aligned}
|< D_\eps^k \tilde U| D_\eps^k \tilde Q(\tilde U)>| &+|< (\eps \d_1)^\ell \tilde U| (\eps \d_1)^\ell \tilde Q(\tilde U)>|
\\
&\leq C\| \tilde U\|_{W^{1,\infty}_\eps} ( \| D_\eps^2 \tilde U \|_{L^2} + \| (\eps \d_1)^4 \tilde U\|_{L^2}+\|\tilde U\|_{L^2})^2 \\
& \leq \frac C\eps \Bigl( \| D_\eps^2 \tilde U \|_{L^2} + \| (\eps \d_1)^4 \tilde U\|_{L^2}+\|\tilde U\|_{L^2} \Bigr)^3.
\end{aligned}
\end{equation}
\begin{Rmk} \label{rk-embedding} Note that because of the bad embedding inequality
$\displaystyle\| \nabla U\|_{L^\infty} \leq \frac1\eps \| U\|_{W_\eps^{4}}  ,$
we lose one power of  $\eps$, which seems not to be optimal considering for instance the fast oscillating functions $\displaystyle x_2\mapsto \exp \left(\frac{ ik_2 x_2}\eps\right).$
A challenging question in order to apply semiclassical methods to nonlinear problems is to determine appropriate functional spaces (in the spririt of~\cite{JMRwigner}) which measures on the one hand the Sobolev regularity of the amplitudes, and on the other hand the oscillation frequency.
\end{Rmk}
We are finally able {\bf to obtain a uniform life span} for the weakly nonlinear system, thus proving result~(\ref{uniformlifespan}) of Theorem~\ref{nonlineartheorem}. Indeed combining the trilinear estimate (\ref{trilinear}) and the commutator estimate (\ref{commutator}), we obtain the following Gronwall  inequality
$${\eps^2} {d\over dt}\left(  \| D_\eps^2 \tilde U \|^2_{L^2} + \| (\eps \d_1)^4 \tilde U\|^2_{L^2}+\|\tilde U\|_{L^2}^2\right)\leq C\eps^2 \left(1+  \| D_\eps^2 \tilde U \|_{L^2} + \| (\eps \d_1)^4 \tilde U\|_{L^2}+\|\tilde U\|_{L^2}\right)^3$$
from which we deduce the uniform a priori estimate
$$\| D_\eps^2 \tilde U \|^2_{L^2} + \| (\eps \d_1)^4 \tilde U\|^2_{L^2}+\|\tilde U\|^2_{L^2}\leq \left(C_0-Ct \right)^{-2}$$
where $C_0$ depends only on the initial data.
Such an estimate shows that the life span of the  solutions to (\ref{NLS}) is at least $T^*=C_0/C$.
  
\subsection{Approximation by the linear dynamics}$ $
In this paragraph we shall prove results~(\ref{iftheorem}) and~(\ref{weakercoupling}) of Theorem~\ref{nonlineartheorem}. 
 The proof of both results relies on standard energy estimates.
 We have
 $$
\eps^2 \d_t (U_\eps -V_\eps) +A(x,\eps D_x) (U_\eps-V_\eps) +\eps^{3+\eta} \tilde S_1(U_\eps) \eps \d_1 U_\eps +\eps^{3+\eta} \tilde S_2(U_\eps) \eps \d_2 U_\eps=0$$

$\bullet$ If $\eta=0$ and $\eps V_\eps \to 0 $ in $L^\infty$, we use the  decomposition
$$ 
\eps^2 \d_t (U_\eps -V_\eps) +A(x,\eps D_x) (U_\eps-V_\eps) +\eps^{3} (\tilde S_j(U_\eps)-\tilde S_j(V_\eps)) \eps\d_j U_\eps +\eps^{3} \tilde S_j(V_\eps)\eps \d_j U_\eps=0$$ 
and obtain the following  $L^2$ estimate
$$\begin{aligned}
\frac{\eps^2}2 {d\over dt} \|U_\eps -V_\eps\|^2_{L^2}       &\leq  3\eps^3 \| \eps \d_j U_\eps \|_{L^\infty} \| U_\eps-V_\eps\|^2_{L^2} + 3\eps^3 \|V_\eps\|_{L^\infty} \| \eps \d_j U_\eps \|_{L^2} \| U_\eps-V_\eps\|_{L^2} \\
&\leq C\eps^2 (\eps  \| \eps \d_j U_\eps \|_{L^\infty} + \| \eps \d_j U_\eps \|_{L^2}^2) \| U_\eps-V_\eps\|^2_{L^2} + C \eps^2 (\eps \|V_\eps\|_{L^\infty} )^2
\end{aligned}$$
from which we conclude by Gronwall's lemma
$$ 
 \|U_\eps -V_\eps\|^2_{L^2}\leq C\int_0^t (\eps \|V_\eps(s)\|_{L^\infty} )^2 \exp C \left( \int_s^t  (\eps  \| \eps \d_j U_\eps \|_{L^\infty} + \| \eps \d_j U_\eps \|_{L^2}^2)  d\sigma\right) ds $$
 on $[0,T^*[$, and that proves result~(\ref{iftheorem}).

$\bullet$ If $\eta>0$, the same arguments show that the life span of the solutions to (\ref{WNLS}) tends to infinity as $\eps \to 0$:
$T_\eps \geq C \eps^{-\eta}\,,$
and that these solutions are uniformly bounded in $W_\eps^{4}$ on any finite time interval.
Furthermore, on any finite time interval $[0,T]$, the previous energy estimate gives
$$
\frac{\eps^2}2 {d\over dt} \|U_\eps -V_\eps\| ^2_{L^2}       \leq  C\eps^{3+\eta}  \| \eps \d_j U_\eps \|_{L^\infty} \| U_\eps\|^2_{L^2} \,,$$
from which we deduce
$$ 
 \|U_\eps -V_\eps\| ^2_{L^2}\leq C\eps^\eta \int_0^t  \eps \| \eps \d_j U_\eps(s)\|_{L^\infty} \| U_\eps(s) \|_{L^2}^2 ds\,.$$
Result~(\ref{weakercoupling}) of Theorem~\ref{nonlineartheorem} is proved.

\appendix
\section{A diagonalization theorem} \label{sectionlemkiditou}
In this appendix we shall state and prove the crucial theorem allowing to
diagonalize semiclassically the matrix of pseudo-differential operators~$A(x_2,\e D)$. 
The construction of Rossby  and Poincar\'e   modes (see Propositions~\ref{therossbymode} and~\ref{thepoincaremode} in Section~\ref{scalar}) are direct corollaries of that theorem.   We refer to Appendix~\ref{semiclassic} for the notation  and results of semi-classical analysis used in this paragraph.
We consider an order function~$g$
 on~$\R^{2d}$,
 and  a symbol~$h(X ; \eps,\tau)$ in~$S_{2d}(g)$, depending polynomially on~$\eps$ and~$\tau$,   where we have defined~$X = (x,\xi) \in \R^{2d}$.  The statement is the following.
   \begin{Thm}\label{lemkiditou}
Let~$g$ be an order function on~$\R^{2d}$, and let~$h(X ; \eps,\tau)$ be a classical symbol in~$S_{2d}(g)$, depending polynomially on~$\eps$ and~$\tau$:   there is an integer~$N_0$  such that
 \begin{equation}\label{polysymbol}
 \forall X \in \R^{2d}, \quad h(X ; \eps,\tau) = \sum_{j, \ell=0}^{N_0}  h^{j,\ell} (X  ) \eps^\ell \tau^j,
 \end{equation}
 where the symbols~$  {h^{j,\ell} (X  )} $   belong  to~$S_{2d}(g)$. 
 Let~$\tilde\tau_\e = \tilde \tau_\e (X )$ be a root of the polynomial~$h(X;\eps,\cdot)$  which    can be written for some~$\nu \in \N$,
$\displaystyle 
 \tilde\tau_\e   (X  ) = \sum_{k=0}^\infty  \e^{\nu + k}  \tau_k (X  ) + O(\eps^\infty), $  {with} $ \tau_0 \neq 0
$
and where~$ \tau_0$ is a symbol.
Finally let~$h_0 (X,\tau)$ be the principal symbol of~$h$, 
satisfying the following assumption:  
\be\label{hyp0}
\exists C>0,\ \forall X \in \R^{2d}, \ \forall \e \in ]0,1[, \vert\partial_\tau h_{0}(X;\e^\nu \tau_0)\vert\geq C.
\ee
Let~${\mathcal K}$ be a compact subset of~$\R^{2d}$. Then there is a pseudo-differential operator $T$ of principal symbol~$ \eps^\nu \tau_0  $ such that if~$\psi \in {\mathcal S}'(\R^d)$ is   microlocalized   in~${\mathcal K}$ and satisfies~$i\partial_s \psi = T\psi $ then
\be
\label{theother}
   {\rm Op}_\eps^0(h (\eps, i \partial_s)) \psi=O(\varepsilon^\infty) \quad \mbox{in} \quad L^2(\R^d).
\ee
\end{Thm}
\begin{proof}
The idea of the proof of the theorem
is the following. Let us define a smooth function~$\chi$, compactly supported in~$\R^{2d}$, identically equal to one on~${\mathcal K}$. Then for any integer~$N$, we shall compute
recursively the coefficients of the symbol
$
   \tau_\e ^N:= \e^\nu \tau_0 + \sum_{k = 1}^N \e^{\nu + k} \tau_k^-  
   $
   so~$T_{\chi} := {\rm Op}_\eps^1 (\tau_{\e}^\chi)  $ satisfies the required property, where~$\tau_\e^\chi$, unique up to~$O(\e^\infty)$, is given by
   $$
   \tau_{\e}^\chi:= \e^\nu \tau_0 + \sum_{k = 1}^\infty \e^{\nu + k}  \chi\tau_k^- + O(\e^\infty) .
   $$
   This allows to replace the root~$\widetilde \tau_\e$ by an actual symbol.
    The above strategy will be achieved in the following way. We notice that   if~$i\partial_s \psi = T_{\chi}\psi $    then of course
    $$
    i\partial_s \psi = {\rm Op}_\eps^1 ( \tau_{\e}^{\chi,N}) \psi + O(\e^{N+\nu+1 }) , \quad \mbox{where} \quad 
   \tau_{\e}^{\chi,N}:= \e^\nu \tau_0 + \sum_{k = 1}^N \e^{\nu + k}  \chi\tau_k^- 
.
 $$
   If moreover~$\psi $   is  microlocalized   in~${\mathcal K}$,  then by definition of~$\chi$  one has   using~(\ref{moyal}-\ref{dieze}) and the fact that~$\chi$ is identically equal to one over~${\mathcal K}$ \   $
\displaystyle i\partial_s \psi =  {\rm Op}_\eps^1 (\tau_{\e}^N  ) \psi + O(\e^{N+\nu+1 }).
 $ 
 This means one can (and shall) compute recursively~$ (\tau_k^- )_{1 \leq k \leq N}$ so that
 $$
    {\rm Op}_\eps^0(h (\eps, i \partial_s)) \psi=O(\varepsilon^{N +1 }) , \quad \mbox{when} \quad 
i\partial_s \psi =  {\rm Op}_\eps^1 (\tau_{\e}^N  ) \psi + O(\e^{N +1 }).
 $$
 Note that it is convenient in the computations to compute~$T_{\chi}$ as the ``right"-quantization   of  the symbol~$ \tau_\e ^{\chi}$. 
Now let us carry out the algebraic computations allowing to achieve the result.
We shall start by dealing with the case when~$\nu \neq 0$, as the computations can be carried out in an easier way, and then we shall discuss the case when~$\nu = 0$.
Recalling that~$h(X;\eps,\tilde\tau_\e(X)) = 0$ we infer that~$h^{0,\ell} \equiv 0$ if~$\ell < \nu $ and
\be\label{order1}
 h^{0,\nu}(x,\xi)  + h^{10} (x,\xi)\tau_0 (x,\xi)= 0. 
\ee
Now we recall that
  for any~$\psi$ in~${\mathcal S}'(\R^d)$, one has
$$
 {\rm Op}_\eps^0(h (\eps,  i\partial_s )) \psi (x)= (2\pi \eps)^{-d} \int_{\R^{2d}} e^{i\frac{(x-y) \cdot \xi}{\eps}}
 h (x,\xi;\eps,  i\partial_s ) \psi(y )\: dyd\xi .
$$
The above integral, as all the ones appearing in this proof, is to be understood in the distributional sense.
So if~$\tau_\e$ and~$\psi$ are such that~$   {\rm Op}_\eps^1 (\tau_\e ^{\chi}) \psi = i\partial_s   \psi$ with~$\psi $ microlocalized   in~${\mathcal K}$,  
\begin{eqnarray}\label{Ope0h} 
{\rm Op}_\eps^0(h (\eps,i\partial_s   )) \psi (x) & = & 
  (2\pi \eps)^{-d} \int_{\R^{2d}} e^{i\frac{(x-y) \cdot \xi}{\eps}}
 h (x,\xi;\eps,  {\rm Op}_\eps^1 (\tau_\e ^\chi)) \psi(y) \: dyd\xi \nonumber \\
 & = & 
  (2\pi \eps)^{-d} \int_{\R^{2d}} e^{i\frac{(x-y) \cdot \xi}{\eps}}
 h (x,\xi;\eps,  {\rm Op}_\eps^1 (\tau_\e^N)) \psi(y) \: dyd\xi + O(\e^{N+1}) .
\end{eqnarray}
 Now we need to compute~$ h (x,\xi;\eps,  {\rm Op}_\eps^1 (\tau_\e^N)) $.  Using the fact that~$h$ is polynomial in~$\tau$,  
 \begin{eqnarray*}
 h (x,\xi;\eps,  {\rm Op}_\eps^1 (\tau_\e^N)) &=& \sum_{j, \ell=0}^{N_0}  h^{j,\ell} (x,\xi ) \eps^\ell \bigl( {\rm Op}_\eps^1 (\tau_\e^N)\bigr)^j \\
 &=:&  {\rm Op}_\eps^1 (\tilde h (x,\xi,y,\eta;\eps))
 \end{eqnarray*}
 where using compositions rules  recalled in Appendix~\ref{semiclassic}, the symbol~$\tilde h (x,\xi,y,\eta;\eps)$ (in the~$(y,\eta)$ variable) can be expanded as
 $$
  \tilde h (x,\xi,y,\eta;\eps) =\e^\nu \sum_{k = 0}^\infty  \tilde h^k (x,\xi,y,\eta) \eps^k + O(\eps^\infty)
 $$
 where the principal symbol of~$  \tilde h (x,\xi,y,\eta;\eps) $ is
 \be\label{pptildeh}
 \e^\nu  \tilde h^{0}(x,\xi,y,\eta)  =  \e^\nu  \left(h^{0,\nu}(x,\xi) + h^{1,0} (x,\xi)\tau_0(y,\eta)   \right).
 \ee
 We recall indeed that according to  Appendix~\ref{semiclassic}, one has
 $
 \displaystyle \bigl( {\rm Op}_\eps^1 (\tau_\e^N)\bigr)^j =:  {\rm Op}_\eps^1 (m_j (y,\eta;\eps)),
 $
where~$\displaystyle m_1 (y,\eta;\eps) = \tau_\e^N (y,\eta)$ and
$$
m_j (y,\eta;\eps) = \sum_{k \geq 0} \frac{(i\eps)^k}{k !} \partial_{y  }^k m_{j-1} (y ,\eta ;\eps)  \partial_{\eta }^k \tau_\e^N (y,\eta) 
+  O(\e^\infty) .
$$
In particular the principal symbol of~$ \bigl( {\rm Op}_\eps^1 (\tau_\e)\bigr)^j$ is~$\e^{\nu j} \tau_0^j(y,\eta)$, which yields~(\ref{pptildeh}).
 We notice that due to~(\ref{order1}), (\ref{pptildeh}) implies that
\be\label{0atxxi}
  \tilde h^{0}(x,\xi,x,\xi) = 0,
 \ee
 so that~$  \tilde h (x,\xi,x,\xi;\eps) = O(\e^{\nu +1})$.
 More generally, plugging the expansion of~$m_j$ into the formula defining~$\tilde h$  and noticing that~$ h^{1,0} (x,\xi) = \partial_\tau h_{0|\tau=0} (x,\xi) $, one finds that there is~$h^k (x,\xi,y,\eta)$, depending on the  symbol coefficients of~$h (x,\xi) $, and in the~$(y,\eta) $ variables on the symbols~$\tau_0 (y,\eta) , \tau^-_{1} (y,\eta), \dots, \tau^-_{k-1}  (y,\eta)$ only, such that
  \be\label{moregenerally}
  \forall k \geq 1,  \quad \tilde h^k (x,\xi,y,\eta) = h^k (x,\xi,y,\eta) +  \partial_\tau h_{0|\tau=0} (x,\xi) \tau^-_k (y,\eta) .
  \ee
  Finally going back to~(\ref{Ope0h}) we find that
 $$
 {\rm Op}_\eps^0(h (\eps, \lambda)) \psi (x)= 
  (2\pi \eps)^{-2d} \int_{\R^{4d}} e^{i\frac{(x-y) \cdot \xi}{\eps}}
  e^{i\frac{(y-y') \cdot \eta}{\eps}}
\tilde h (x,\xi,y',\eta;\eps ) \psi(y') \: dy dy'd\xi d\eta + O(\e^N) 
 .
 $$
 We can first perform the integration in the~$y$ variable, which creates a Dirac mass at~$\xi -\eta$, and therefore we have
  \be\label{diracy}
 {\rm Op}_\eps^0(h (\eps, \lambda)) \psi (x)= 
  (2\pi \eps)^{-d} \int_{\R^{2d}} e^{i\frac{(x-y') \cdot \xi}{\eps}}
  \tilde h (x,\xi,y',\xi;\eps ) \psi(y') \:  dy'd\xi + O(\e^N) 
  .
 \ee
Now we shall construct~$\tau^-_1$  ensuring that the order of~$ {\rm Op}_\eps^0(h (\eps, \lambda)) \psi (x)$ is~$ O(\e^{\nu +1})$ instead of~$ O(\e^{\nu })$. The argument will easily be adaptable by induction to~$ {\rm Op}_\eps^0(h (\eps, \lambda)) \psi (x)= O(\e^N)$ for any~$N \in \N$, by a convenient choice of~$\tau^-_j$, for~$j \leq N-1$.  We notice that in~(\ref{diracy}), the quantity~$ \tilde h (x,\xi,y',\xi;\eps ) $ can easily be replaced by~$ \tilde h (x,\xi,x,\xi;\eps ) $ 
 by Taylor's formula: more precisely we write
 \be\label{firstorderexpansion}
  \tilde h (x,\xi,y',\xi;\eps ) =  \tilde h (x,\xi,x,\xi;\eps ) + (y'-x) \cdot( \nabla_{y'}  \tilde h) (x,\xi,x,\xi;\eps ) + O(|y'-x|^2)
 \ee
 which gives after integrations by parts
   \begin{eqnarray*}
 {\rm Op}_\eps^0(h (\eps, \lambda)) \psi (x) & = & \int_{\R^{2d}} e^{i\frac{(x-y') \cdot \xi}{\eps}}
  \tilde h (x,\xi,x,\xi;\eps ) \psi(y') \:  \frac{dy'd\xi}{ (2\pi \eps)^{ d} } \\
  && {}  - i\e  \int_{\R^{2d}} \nabla_\xi e^{i\frac{(x-y') \cdot \xi}{\eps}}
   \cdot (\nabla_{y'}  \tilde h) (x,\xi,x,\xi;\eps ) \psi(y') \:  \frac{dy'd\xi}{ (2\pi \eps)^{ d} }    + O(\e^{\nu +2}) \\
   & = &  \int_{\R^{2d}} e^{i\frac{(x-y') \cdot \xi}{\eps}}
  \tilde h (x,\xi,x,\xi;\eps ) \psi(y') \:  \frac{dy'd\xi}{ (2\pi \eps)^{ d} } \\
  && {}  + i\e  \int_{\R^{2d}}e^{i\frac{(x-y') \cdot \xi}{\eps}}
   \nabla_\xi  \cdot(\nabla_{y'}  \tilde h) (x,\xi,x,\xi;\eps ) \psi(y') \:  \frac{dy'd\xi}{ (2\pi \eps)^{ d} }   + O(\e^{\nu +2}).
 \end{eqnarray*}
 Now using~(\ref{0atxxi}), it remains to choose~$\tau^-_1$ so that
    $$
  \int_{\R^{2d}} e^{i\frac{(x-y') \cdot \xi}{\eps}}
 \bigl( \tilde h (x,\xi,x,\xi;\eps )  + i\e     \nabla_\xi  \cdot( \nabla_{y'}  \tilde h) (x,\xi,x,\xi;\eps ) \bigr) \psi(y') \:   \frac{dy'd\xi}{ (2\pi \eps)^{ d} }  = O(\e^{\nu +2}).
   $$
   That is possible
 simply by looking at formula~(\ref{moregenerally})  and choosing
 $$
\tau^-_1 (x,\xi)  = -\frac{ \tilde h^1(x,\xi,x,\xi)  +i  \nabla_\xi  \cdot (\nabla_{y}  \tilde  h^0 )(x,\xi,x,\xi  )}{ \partial_\tau h_{0|\tau=0} (x,\xi) }  \cdotp
$$
Note that Assumption~(\ref{hyp0}) guarantees that~$\tau^-_1 $ is well defined.
  The argument may be     pursued at the next order simply replacing~(\ref{firstorderexpansion}) by
     \begin{eqnarray*}
  \tilde h (x,\xi,y',\xi;\eps ) =  \tilde h (x,\xi,x,\xi;\eps ) &+& (y'-x) \cdot( \nabla_{y'}  \tilde h) (x,\xi,x,\xi;\eps )\\&& {}+ \frac12
  (y'-x) \otimes  (y'-x) : (\nabla_{y'}^2 \tilde h) (x,\xi,x,\xi;\eps ) +
  O(|y'-x|^3)
 \end{eqnarray*}
and using integrations by parts again.  Then the choice
  $$
 \tau^-_2 (x,\xi)  = -\frac{ \tilde h^2(x,\xi,x,\xi)  +i  \nabla_\xi  \cdot( \nabla_{y}   \tilde h^1) (x,\xi,x,\xi ) -  \frac12\nabla_\xi^2 :(\nabla_{y}^2  \tilde  h^0) (x,\xi,x,\xi  ) }{ \partial_\tau h_{0|\tau=0} (x,\xi) } 
 $$
gives that~$  {\rm Op}_\eps^0(h (\eps, \lambda)) \psi (x) = O(\e^{\nu +3})$. We leave the rest of the induction argument to the reader.
To end the proof of Theorem~\ref{lemkiditou}
 we   need to consider the case~$\nu = 0$. The argument is similar to the case~$\nu \neq 0$ treated above, though the formulas are slightly more complicated. We shall use the same notation as in the previous case.
We start by noticing that  by definition of~$\tilde \tau_\e$ we have in particular  for all
$
 (x,\xi) \in \R^{2d}, $ $ \displaystyle  \sum_{j = 0}^{N_0} h^{j,0} (x,\xi) \tau_0^j (x,\xi)= 0.
$
Then  as before let us write
$
\displaystyle  h (x,\xi;\eps,  {\rm Op}_\eps^1 (\tau_\e^N))  =:   {\rm Op}_\eps^1 (\tilde h (x,\xi,y,\eta;\eps)).
$
One computes easily that   the principal symbol of~$  \tilde h (x,\xi,y,\eta;\eps) $ is (unlike the case~$\nu \neq 0$)
 \be\label{pptildeh0}
   \tilde h^{0}(x,\xi,y,\eta)  =  \sum_{j=0}^{N_0} h^{j,0}(x,\xi) \tau_0^j (y,\eta)  .
 \ee
 Note that as above one has
 $
  \tilde h^{0}(x,\xi,x,\xi) = 0,
 $
  so that~$  \tilde h (x,\xi,x,\xi;\eps) = O(\e)$.
  One can also compute the next orders, and   as in the case~$\nu \neq 0$ they can be written in the following form:
  $$
   \forall k \geq 1,  \quad \tilde h^k (x,\xi,y,\eta) = h^k (x,\xi,y,\eta) +  \tau^-_k (y,\eta)  \sum_{j=1}^{N_0} h^{j,0}(x,\xi) j \tau_0^{j-1} (y,\eta) 
 $$
   where~$h^k (x,\xi,y,\eta)$ depends on the  symbol coefficients of~$h (x,\xi) $, and in the~$(y,\eta) $ variables on~$\tau_0 (y,\eta) , \tau^-_{1} (y,\eta), \dots, \tau^-_{k-1}  (y,\eta)$ only. In particular we notice that
      \be\label{moregenerally0}
\forall k \geq 1,  \quad \tilde h^k (x,\xi,x,\xi) = h^k (x,\xi,x,\xi) +    \tau^-_k (y,\eta) \partial_\tau h_{0|\tau = \tau_0}.
 \ee
Now that these formulas have been established, it remains to go through exactly the same computations as in the case~$\nu = 0$, and we find
$$
\tau^-_1 (x,\xi)  = -\frac{\tilde h^1(x,\xi,x,\xi)  +i  \nabla_\xi  \cdot (\nabla_{y}  \tilde  h_0 )(x,\xi,x,\xi  )}{ \partial_\tau h_{0|\tau=\tau_0} (x,\xi) }  $$
which is well defined thanks to Assumption~(\ref{hyp0}). The other orders are obtained exactly as in the case~$\nu = 0$. 
This ends the proof of Theorem~\ref{lemkiditou}.
   \end{proof}

\section{Some well-known facts in semi-classical analysis}\label{semiclassic}
In this section we recollect some well-known facts in semi-classical analysis, which have been used 
throughout the paper. Most of the material is taken from~\cite{DS}, \cite{hormander}, \cite{MA}, \cite{san} 
and~\cite{stein}.
\subsection{Semi-classical symbols and operators}
\subsubsection{Definitions}
We   recall that an {\bf order function} is any function~$g \in C^\infty(\R^{d};\R^+\setminus\{0\})$ such that 
there is a constant~$C$ satisfying
$$
\forall X \in \R^d, \:  \forall  \alpha \in \N^d, \quad |\partial^\alpha g(X)| \leq C g(X).
$$
For instance~$ g(x,\xi) = (1+|\xi|^2)^\frac12 =: \langle \xi \rangle$ is an order function.
Note that the variable~$X$ usually refers  to a point~$(x,\xi)$ in the cotangent space~$T^*\R^n \equiv \R^{2n}$, or to a point of the type~$(x,y,\xi)$ with~$y \in \R^n$. 
A {\bf semi-classical symbol} in the class~$S_d(g)$ is   then a function~$a = a(X;\varepsilon)$ defined on~$\R^{d} \times ]0,\eps_0]$ for some~$\eps_0 >0$, which depends smoothly on~$X$ and such that for any~$\alpha \in \N^d$, there is a constant~$C$ such that 
$  |\partial^\alpha a(X,\eps)| \leq C g(X) $ for any~$(X,\eps) \in \R^{d} \times ]0,\eps_0]$.

If~ $(a_j)_{j \in \N}$ is a family of  semi-classical symbols in the class~$S_d(g)$, we write that
$$
a = \sum_{j=0}^\infty \eps^j a_j + O(\eps^\infty)
$$
if for any~$N \in \N$ and for any~$\alpha \in \N^d$, there are~$\e_0$ and~$C$ such that
$$
\forall  X \in \R^d, \: \forall \eps \in  ]0,\eps_0]\quad  \Bigl|\partial^\alpha \Bigl(a(X,\eps) - \sum_{j=0}^N \eps^j a_j(X,\eps) \Bigr) \Bigr| \leq C \e^Ng(X).
$$
Conversely for any sequence~$(a_j)_{j \in \N}$ of symbols in~$S_d(g)$, there is~$a \in S_d(g)$ (unique up to~$ O(\eps^\infty)$) such that
$
\displaystyle a= \sum_{j=0}^\infty \eps^j a_j + O(\eps^\infty).
$
An {\bf $\eps$-pseudodifferential operator}  is      defined as follows: if~$a$ belongs to~$S_{3n}(g)$, and~$u$ is in~$ \mathcal D(\R^n)$, then
$$
\Bigl({\rm Op}_\eps(a) \Bigr)u(x ) := \frac1{(2\pi \eps)^n} \int e^{i(x-y) \cdot \xi/\eps} a(x,y,\xi) u(y) \: dyd\xi
.$$
\subsubsection{Changes of quantization}
If~$a \in S_{2n}(g)$ and~$t \in [0,1]$ then~$
a^t(x,y,\xi) := a( (1-t) x + ty,\xi)
$
belongs to~$S_{3n}(g)$, and one defines $ {\rm Op}^t_\eps(a) := {\rm Op}_\eps (a^t) $.
When~$t = 0$ this corresponds to the {\bf classical}, or ``{\bf left}" quantization, when~$t = 1/2$ this is known as the {\bf Weyl quantization} (and is usually denoted by~${\rm Op}^W_\eps(a) = {\rm Op}^\frac12_\eps(a) $), while when~$t=1$ one refers to the ``{\bf right}" quantization. Furthermore, if~$a$ belongs to~$S_{3n}(\langle \xi\rangle^m)$ for some integer~$m$ (or more generally if~$a$ belongs to~$S_{3n}(g)$ where~$g$ is a H\"ormander metric~\cite{hormander}), then there is a unique symbol~$a_t $ belonging to~$ S_{2n}(\langle \xi\rangle^m)$ (resp. $a_t \in S_{3n}(g)$) such that~${\rm Op}^t_\eps(a_t)  = {\rm Op}_\eps (a)$, and one has
$$
a_t(x,\xi) = \sum_\alpha \frac{(i\eps)^\alpha}{\alpha !} \partial_\xi^\alpha\partial_\eta^\alpha
a(x+t\eta,x-(1-t)\eta,\xi)_{|\eta = 0} + O(\eps^\infty).
$$
A {\bf classical symbol} is a symbol~$a$ in~$S_{2n}(\langle \xi\rangle^m)$ such that
$
\displaystyle a(x,\xi ; \eps) = \sum_{j=0}^\infty \eps^j a_j(x,\xi) + O(\eps^\infty)
$
with~$a_0$ not identically zero, and~$a_j \in S_{2n}(\langle \xi\rangle^m)$ independent of~$\eps$. For any real number~$\nu$, the term~$\eps^\nu a_0$ is the {\bf principal symbol} of the classical pseudo-differential operator~$A = \eps^\nu {\rm Op}^t_\eps(a) $ (and this does not depend on the quantization).  On the other hand~$\eps^{\nu+1} a_1$ is the   {\bf subprincipal} symbol of~$A = \eps^\nu {\rm Op}^W_\eps(a) $ (in the Weyl quantization only).
In the following we shall denote by~$\sigma_t (A)   $   the symbol of an operator~$A =  {\rm Op}^t_\eps(a) $ (in other words~$a = \sigma_t (A) $), and by~$\sigma_P(A)$ its principal symbol.
\subsubsection{Microlocal support and~$\e$-oscillation}
If~$u$ is an~$\eps$-dependent function in a ball of~$L^2(\R^n)$, its~{\bf $\eps$-frequency set} (or {\bf  microlocal support}) is the complement in~$\R^{2n}$ of the points~$(x_0,\xi_0)  $ such that there is a function~$\chi_0 \in S_{2n}(1)$ equal to one at~$(x_0,\xi_0)$, satisfying
$$
\| {\rm Op}^W_\eps (\chi_0u) \|_{L^2(\R^{2n})} = O(\eps^\infty).
$$
We   say that an~$\eps$-dependent function $f_\e$ bounded in~$L^2(\R^n)$ is~{\bf $\e$-oscillatory} if for every continuous, compactly supported function~$\varphi$ on~$\R^n$, 
\begin{equation}\label{epsosc}
\limsup_{\e \to 0} \int_{|\xi| \geq R/\e} |\varphi \widehat f_\e  (\xi)|^2 \: d\xi \to 0 \quad \mbox{as} \quad R \to \infty.
\end{equation}
An~$\eps$-dependent function $f_\e$ bounded in~$L^2(\R^n)$ is said to be {\bf compact at infinity} if
\begin{equation}\label{compinfinity}
\limsup_{\e \to 0} \int_{x \geq R } | f_\e(x)|^2 \: dx \to 0 \quad \mbox{as} \quad R \to \infty.
\end{equation}
 \subsubsection{Adjoint and composition}
 Let $a$  be a symbol  in~$S_{3n} (g) $, where~$g$ is a H\"ormander metric~\cite{hormander}, and define
 $
 a^* (x,y,\xi) := \overline{a  (y,x,\xi) }.
 $
 Then the operator~$( {\rm Op}_\eps (a) )^* :=   {\rm Op}_\eps (a^*)$ satisfies for all~$u,v$ in~${\mathcal S}(\R^n)$, 
 $$
 \Bigl( ( {\rm Op}_\eps (a) )^* u , v \Bigr)_{L^2} =  \Bigl( u ,  ( {\rm Op}_\eps (a) )^*v \Bigr)_{L^2} 
 $$
 and is therefore called the {\bf formal adjoint} of~$  {\rm Op}_\eps (a)  .$ In particular~$  {\rm Op}^\frac12_\eps (a) $ is  {\bf formally self-adjoint} if~$a$ is real.  
Let $a$ and~$b$ be two symbols in~$S_{2n} (g_1) $ and~$S_{2n} (g_2) $ respectively, where~$g_j$ are H\"ormander metrics. For all~$t \in [0,1]$, there is a unique symbol~$c_t$ in~$S_{2n} (g_1g_2 ) $ which allows to obtain 
 $ {\rm Op}^t_\eps(a) \circ {\rm Op}^t_\eps(b) = {\rm Op}^t_\eps(c_t) $. Moreover one has
 \begin{equation}\label{moyal}
 c_t(x,\xi;\e) = e^{i\e [ \partial_u  \partial_\xi - \partial_\eta  \partial_v ]} \left(
 a((1-t)x + tu,\eta) b((1-t)v + tx,\xi)
 \right)_{ {\scriptsize \begin{array}{c}
  | u=v=x\\ 
 \eta = \xi 
 \end{array}}
 } =: a \#^t b.
 \end{equation}
  This can be also written
   $$
 c_t(x,\xi;\e) =\displaystyle  \sum_{k \geq 0} \frac{\e^k}{i^k k !} (\partial_\eta \partial_v - \partial_\xi \partial_u)^k \left(
 a((1-t)x + tu,\eta) b((1-t)v + tx,\xi)
 \right)_{ {\scriptsize \begin{array}{c}
  | u=v=x\\ 
 \eta = \xi 
 \end{array}}
 }  + O(\eps^\infty) .
 $$
In particular one has
$
\sigma_t (A \circ B)  = 
\sigma_t (A) \sigma_t (B) + O(\e).
$
  For example in the case when~$t=0$ then~$\displaystyle {\rm Op}_\eps(a) \circ {\rm Op}_\eps(b) = {\rm Op}_\eps(c) $, with
   \begin{equation}\label{dieze}
  c(x,\xi) =  a \#  b =  \sum_{\alpha} \frac{\e^ \alpha}{i^ \alpha \alpha !}  \partial_\xi^\alpha a(x,\xi)  \partial_x^\alpha b(x,\xi)  + O(\eps^\infty) .
 \end{equation}
  In particular if~$a$ and~$b$ are two classical symbols in the sense described above, and if one defines~$A :=  {\rm Op}^t_\eps(a)$ and~$B=  {\rm Op}^t_\eps(b)$, then the principal symbols satisfy (if~$\sigma_P(A)\sigma_P(B $ does not vanish identically)
$
 \sigma_P(AB) = \sigma_P(A)\sigma_P(B).
$
  \subsection{Semiclassical operators, Wigner transforms  and   propagation of   energy}
One of the main interests of the semiclassical setting is that it allows a precise description of the propagation of the energy, on times of the order of~$O(\e)$. We refer for instance to~\cite{gmmp} (Section 6) for the proof of the following property (actually in the more general setting of matrix-valued operators): consider a scalar symbol~$\tau_\e(x,\xi)$  defined on~$\R^{2n}$, belonging to the class~$S_{2n} (\langle \xi\rangle^\sigma) $ for some~$\sigma \in \R$ (or more generally to~$S_{2n} (g) $ where~$g$ is a H\"ormander metric). We assume moreover that~$ {\rm Op}_\e^W (\tau_\e)$ is essentially skew-self-adjoint on~$L^2(\R^n)$.  Then consider~$f_\e^0$   an~$\e$-oscillatory initial data in the sense of~(\ref{epsosc}),  
bounded in~$L^2(\R^n)$ and compact at infinity in the sense of~(\ref{compinfinity}), and the PDE
$$
\e \partial_t f_\e + {\rm Op}_\e^W (\tau_\e) f_\e = 0, \quad f_{\e|t = 0} = f_\e^0.
$$
Then the Wigner transform~$W_\e(t,x,\xi)$ of~$f_\e(t)$ defined by
$$
W_\e(t,x,\xi) := (2\pi)^{-n} \int_{\R^n} e^{i v \cdot \xi}f_\e(x-\frac \e 2  v) \bar f_\e(x+\frac \e 2  v) \: dv
$$  converges, locally uniformly in~$t$, to the solution~$W$ of
$\partial_t W + \{ \tau_0, W\} =0
$
where~$\tau_0$ is the principal symbol of~$\tau_\e$, and where the Poisson bracket is given by
$$
\{ \tau_0, W\} := \nabla_\xi \tau_0 \cdot \nabla_x W -  \nabla_x \tau_0 \cdot \nabla_\xi W.
$$
The interest of Wigner transforms lies in particular in the fact that under the assumptions made on~$f_\e^0$, for any compact set~$K \subset \R^n$ one has $ \int_{K} |f_\e(t,x)|^2 \: dx = W_\e(t,K \times \R^n) $ due to the fact that~$  |f_\e(t,x)|^2 = \int_{\R^n}W_\e(t,x,\xi  ) \: d\xi$.

\subsection{Coherent states}
A {\bf coherent state}   is  $\displaystyle
\Phi_{p,q}(y) := (\pi\eps)^{-\frac{n}4} e^{i\frac{(y-q) \cdot p }{\eps} } e^{-\frac{(y-q)^2}{2\eps}}
$. Any tempered distribution~$u$ defined on~$\R^n$ may be written
$$
u (y)=   (2\pi \e) ^{-\frac{n}2}  \int Tu(p,q) \Phi_{p,q}(y) \: dpdq,
$$
where~$T$ is the~$FBI$ (for Fourier-Bros-Iagolnitzer) transform
$$
Tu(p,q) :=  2^{-\frac{n}2}  {(\pi \eps)}^{-\frac{3n}4} \int  e^{i\frac{(q-y) \cdot p }{\eps} }  e^{-\frac{(y-q)^2}{2\eps}} u(y) \: dy.
$$
This transformation
 maps isometrically~$L^2(\R^n)$ to~$L^2(\R^{2n})$. 
The above formula simply translates the fact that~$u = T^*Tu$.
 \subsection{Fourier Integral Operators}
A Fourier Integral Operator (FIO) is an operator $U$ which can be written, for any~$f \in L^2(\R^n)$
$$
Uf (x) =\frac1{(2\pi \e)^\frac{3n}{2}} \int_{\R^{2n}} e^{i \Phi (x,y,\tau)/\e} a(x,y,\tau)   f(y) \: dy d\tau
$$
where~$a$ is a   symbol of order~0, compactly supported in~$x$ and~$y$, $\Phi$ is real valued and homogeneous of degree 1 in $\tau$, smooth for~ $\tau \neq 0$. One requires also a non degeneracy condition on the phase (see~\cite{stein}, Chap. 9, Par. 6.11) 
 on the support of~$a$. Then~$U$ is continuous over~$ L^2(\R^n)$.

{\bf Acknowledgements.}  $ $ The authors are grateful to B. Texier for having answered many questions concerning geometric optics and semi-classical analysis. They also thank P. G\'erard for interesting discussions around the  study of Poincar\'e waves in Section~\ref{sectionpoincare}. Finally they thank the anonymous referees for a very careful reading of the manuscript.
I. Gallagher and L. Saint-Raymond are 
partially supported by the French Ministry of Research grant
ANR-08-BLAN-0301-01.

\end{document}